\documentclass[12pt,leqno]{amsart}
\usepackage{amsmath,amsthm,amsfonts,amsbsy,multicol}

\usepackage[dvips]{graphicx}
\input{epsf.tex}
\usepackage[dvipsnames,usenames]{color}
\DeclareGraphicsExtensions{.pdf, .jpg, .pnp, .bmp}

\usepackage{amssymb}
\usepackage{latexsym}

\newtheorem{lemma}{Lemma}[section]

\newtheorem{theorem}[lemma]{Theorem}
\newtheorem{proposition}[lemma]{Proposition}
\newtheorem{definition}[lemma]{Definition}
\newtheorem{corollary}[lemma]{Corollary}
\newtheorem{example}[lemma]{Example}
\newtheorem{exercise}[lemma]{Exercise}
\newtheorem{remark}[lemma]{Remark}
\newtheorem{fig}[lemma]{Figure}
\newtheorem{tab}[lemma]{Table}

\newcommand{\bth}{\begin{theorem}}   \newcommand{\ethe}{\end{theorem}}
\newcommand{\bre}{\begin{remark}\em }   \newcommand{\ere}{\end{remark}}
\newcommand{\ble}{\begin{lemma}}      \newcommand{\ele}{\end{lemma}}
\newcommand{\bde}{\begin{definition}}   \newcommand{\ede}{\end{definition}}
\newcommand{\bco}{\begin{corollary}}     \newcommand{\eco}{\end{corollary}}
\newcommand{\bpr}{\begin{proposition}}  \newcommand{\epr}{\end{proposition}}
\newcommand{\bexer}{\begin{exercise}}     \newcommand{\eexer}{\end{exercise}}
\newcommand{\bexam}{\begin{example}}    \newcommand{\eexam}{\end{example}}
\newcommand{\bfi}{\begin{fig}}              \newcommand{\efi}{\end{fig}}
\newcommand{\btab}{\begin{tab}}         \newcommand{\etab}{\end{tab}}
\newcommand{\bpf}{\begin{proof}}        \newcommand{\epf}{\end{proof}}

\newcommand{\barr}{\begin{array}}   \newcommand{\earr}{\end{array}}
\newcommand{\beao}{\begin{eqnarray*}}       \newcommand{\eeao}{\end{eqnarray*}\noindent}
\newcommand{\beam}{\begin{eqnarray}}    \newcommand{\eeam}{\end{eqnarray}\noindent}
\newcommand{\beqq}{\begin{equation}}    \newcommand{\eeqq}{\end{equation}\noindent}

\newcommand{\bce}{\begin{center}}   \newcommand{\ece}{\end{center}}

\newcommand{\ov}{\overline} 

\newcommand{\wh}{\widehat}
     
  \newcommand{\Nto}{N\to\infty}

\newcommand{\Rto}{R\to\infty}

\newcommand{\al}{\alpha}
\newcommand{\be}{\beta}
 
\newcommand{\del}{\delta}
\newcommand{\D}{\Delta}
  \newcommand{\ep}{\epsilon}
\newcommand{\ka}{\kappa}
\newcommand{\lam}{\lambda} 

\newcommand{\w}{\omega} \newcommand{\W}{\Omega}

\newcommand{\bbc}{{\mathcal C}} 

\newcommand{\bfE}{{\mathbb E}} 
\newcommand{\bbf}{{\mathcal F}}

\newcommand{\bbi}{{\mathbb I}}

\newcommand{\bbl}{{\mathcal L}}
\newcommand{\bbm}{{\mathcal M}}
 \newcommand{\bbN}{{\mathbb N}}
\newcommand{\bfP}{{\mathbb P}}
                              
 \newcommand{\bbR}{{\mathbb R}}

\begin{document}

\bibliographystyle{plain}
\title[On the numerical solution of  super linear stochastic differential equations]{On the numerical solution of some nonlinear stochastic differential equations using the semi-discrete method}

\author[N. Halidias]{N. Halidias}
\author[I. S. Stamatiou]{I. S. Stamatiou}
\address{Department of Mathematics,
University of the Aegean\\ Karlovassi, GR-83 200 Samos, Greece\\ Tel.  +3022730-82321, +3022730-82343}
\email{nick@aegean.gr, istamatiou@aegean.gr}

\begin{abstract}
In this paper we are interested in the numerical solution of
stochastic differential equations with non negative solutions. Our
goal is to construct explicit numerical schemes that preserve
positivity, even for super linear stochastic differential equations.
It is well known that the usual Euler scheme diverges on super linear problems and
the Tamed-Euler method does not preserve positivity.  In that
direction, we use the Semi-Discrete method that the first author has proposed in
two previous papers. We propose a new
numerical scheme for a class of stochastic differential equations
which are super linear with non negative solution. In this class of stochastic differential equations belongs the Heston $3/2$-model that appears in financial
mathematics, for which we prove 
through numerical experiments the ``optimal'' order of strong convergence at least $1/2$ of the Semi-Discrete method.
\end{abstract}
\date\today

 \maketitle

\tableofcontents
\listoffigures
\listoftables

\section{Introduction.}\label{sec:s1}
\setcounter{equation}{0}

Throughout, let $T>0$ and $(\Omega, \bbf, \{\bbf_t\}_{0\leq t\leq T}, \bfP)$ be a complete probability space, meaning that the filtration $ \{\bbf_t\}_{0\leq t\leq T} $ satisfies the usual conditions, i.e. is right continuous and $\bbf_0$ includes all $\bfP-$null sets. Let $W_{t,\w}:[0,T]\times\W\rightarrow\bbR$ be a one dimensional Wiener process adapted to the filtration $\{\bbf_t\}_{0\leq t\leq  T}.$ Consider the following stochastic differential equation (SDE),
\beqq\label{eq1}
x_t=x_0 + \int_{0}^{t}a(s,x_s)ds + \int_{0}^{t}b(s,x_s)dW_s,\quad t\in [0,T],
\eeqq
where the coefficients $a,b: [0,T]\times \bbR\mapsto\bbR$ are measurable functions such that (\ref{eq1}) has a unique strong solution and $x_0$ is independent of all $\{W_t\}_{0\leq t\leq T}$, $x_0>0,$ a.s. SDE (\ref{eq1}) has non autonomous coefficients, i.e. $a(t,x), b(t,x)$ depend explicitly on $t.$

To be more precise, we assume the existence of a predictable stochastic process $x:[0,T]\times \W\mapsto \bbR$ such that (\cite[Definition 2.1]{mao:1997}),
$$
\{a(t,x_t)\}\in\bbl^1([0,T];\bbR), \quad \{b(t,x_t)\}\in\bbl^2([0,T];\bbR)
$$
and
$$
\bfP\left[x_t=x_0 + \int_{0}^{t}a(s,x_s)ds + \int_{0}^{t}b(s,x_s)dW_s\right]=1, \quad \hbox{ for every } t\in[0,T].
$$
The drift coefficient $a$ is the infinitesimal mean of the process $x_t$ and the diffusion coefficient $b$ is the infinitesimal variance of the process $x_t.$ SDEs of the form (\ref{eq1}) have rarely explicit solutions, thus numerical approximations are necessary for simulations of the paths $x_t(\w),$ or for approximation of functionals of the form $\bfE F(x),$ where $F:\bbc([0,T],\bbR)\mapsto\bbR$ can be for example in the area of finance, the discounted payoff of European type derivative.

We are interested in strong approximations (mean-square) of (\ref{eq1}), in the case of super or sub linear drift and diffusion coefficients. This kind of numerical schemes have applications in many areas, such as simulating scenarios, filtering, visualizing stochastic dynamics (see for instance \cite[Section 4]{hutzenhaler_jentzen:2012} and references therein), have theoretical interest (they provide fundamental insight for weak-sense schemes) and generally do not involve simulations over long-time periods or of a significant number of trajectories.

We present some models that are not linear both in the drift and diffusion coefficient:

\begin{itemize}

  \item  The following linear drift model had been initially proposed for the dynamics of the inflation rate in (\cite[Relation 50]{cox_et_al:1985}) and has taken its name, CIR, by the initials of the authors in the aforementioned paper. It is used in the field of finance as a description of the stochastic volatility procedure  in the Heston model (\cite{heston:1993}), but also belongs to the fundamental family of SDEs that approximate  Markov jump processes (\cite{either_kurtz:1986}).  The CIR model is described by the following SDE,
\beqq\label{eq1004}
x_t=x_0 + \int_{0}^{t}\ka(\lam-x_s)ds + \int_{0}^{t}\sigma\sqrt{x_s}dW_s, \quad t\in [0,T],
\eeqq
where $x_0$ is independent of all $\{W_t\}_{t\geq0}, x_0>0,$ a.s. and the parameters $\ka, \lam, \sigma$ are positive. Parameter $\lam$ is the level of the interest rate $x_t$ where the drift is zero, meaning that when $x_t$ is below $\lam$ the drift is positive, whereas in the other case is negative. As $\lam$ grows, the range of the positive drift becomes wider. Parameter $\ka$ defines the slope of the drift. The condition $\ka>0$ is necessary for the stationarity of the process $x_t.$ When $\ka$ is negative, the main term of the slope, $-\ka,$ is positive and  given the diffusion $\sigma\sqrt{x_t},$ the process $x_t$ blows up. The condition $\sigma^2< 2\ka\lam$ implied by the Feller test (\cite[Case (ii),p.173]{feller:1951}) is necessary and sufficient for the process not to reach the boundary zero in finite time.

 \item The $3/2-$model (\cite{heston:1997}) or the inverse square root process (\cite{ahn_gao:1999}), that is used for modeling stochastic volatility,
\beqq\label{eq1003}
x_t=x_0 + \int_{0}^{t}(\al x_s -\be x_s^{2})ds + \int_{0}^{t}\sigma x_s^{3/2} dW_s,\quad t\in [0,T],
\eeqq
where $x_0$ is independent of $\{W_t\}_{0\leq t\leq T}$, $x_0>0,$ a.s. and $ \sigma\in\bbR$. The conditions
 $\al>0$ and $\be>0$ are necessary and sufficient for the stationarity of the process $x_t$ and such that zero and infinity is not attainable in finite time (\cite[Appendix A]{ahn_gao:1999}).

\item The constant elasticity of variance model (\cite{cox:1975}), which is used for pricing assets,
\beqq\label{eq1008}
x_t=x_0 + \int_{0}^{t}\mu x_s ds + \int_{0}^{t}\sigma x_s^{\gamma} dW_s,\quad t\in [0,T],
\eeqq
where $x_0$ is independent of $\{W_t\}_{0\leq t\leq T}$, $x_0>0,$ a.s., $\mu\in\bbR, \sigma>0$ and $0<\gamma\leq1$. SDE (\ref{eq1008}) has a unique strong solution if and only if $\gamma\in[1/2,1]$ and takes values in $[0,\infty).$ The case $\gamma=1/2$ corresponds to CIR model (\ref{eq1004}), whereas $\gamma=1$ corresponds to a Brownian motion, i.e. the famous Black-Scholes model (\cite{black_scholes:1973}).

\item  Superlinear models, i.e. models of the form (\ref{eq1}) where one of the coefficients $a(\cdot),b(\cdot)$  is superlinear, i.e. when we have that
\beqq\label{eq1005}
a(x)\geq \frac{|x|^{\beta}}{C},\, b(x)\leq C|x|^{\alpha}, \, \hbox{for every}\quad |x|\geq C,
\eeqq
or
\beqq\label{eq1006}
b(x)\geq \frac{|x|^{\beta}}{C}, \,a(x)\leq C|x|^{\alpha}, \, \hbox{for every}\quad |x|\geq C,
\eeqq
where $\beta >1, \beta>\alpha\geq0, C>0.$
\end{itemize}

For some of the aforementioned problems there are methods of simulation (\cite{broadie_kaya:2006},  \cite{marakov_glew:2010}). However, if a full sample path of the SDE has to be simulated or the SDEs under study are a part of a bigger system of SDEs, then numerical schemes are in general more effective.

Problems like (\ref{eq1004}), (\ref{eq1003}) and (\ref{eq1008}) are meant for non-negative values, since they represent rates or pricing values. Thus ``good'' numerical schemes preserve positivity (\cite{appleby_et_al:2010}, \cite{kahl_et_al:2008}). The explicit Euler scheme has not that property, since its increments are conditional Gaussian. For example, the transition probability of the Euler scheme in case of (\ref{eq1004}) reads as
$$
p(y|x)=\frac{1}{\sqrt{2\pi\sigma^2x\D}}\exp\Big\{-\frac{(y-(x+\ka(\lam-x)\D))^2}{2\sigma^2x\D}\Big\}, \qquad y\in\bbR,\, x>0,
$$
thus, even in the first step there is an event of negative values with positive probability. We refer to  (\cite{kloeden_neuenkirch:2012}), between other papers, that considers Euler type schemes, modifications of them to overcome the above drawback, and the importance of positivity.  Thus, for the same problem, the truncated Euler scheme (\cite{deelstra_delbaen:1998}) has been proposed, as well as a modification of it, (\cite{higham_mao:2005}), where in a step the numerical scheme can leave $(0,\infty)$ but is forced to come back in the next steps.

One more drawback, that appears in case of superlinear problems (\ref{eq1005}) or (\ref{eq1006}), like (\ref{eq1003}), is that the moments of the scheme may explode (\cite[Theorem 1]{hutzenhaler_et_al.:2011}). A method that overcomes this drawback is the Tamed-Euler method,  (\cite[Relation 4]{hutzenhaler_jentzen:2012}) and reads:  $Y_{0}^N(\w):=x_0(\w)$ and
\beqq\label{eq1000}
Y_{n+1}^N(\w):= Y_n^N(\w) + \frac{T/N \cdot a(Y_n^N(\w)) + b(Y_n^N(\w))\left(W_{\frac{(n+1)T}{N}}(\w) - W_{\frac{nT}{N}}(\w)\right)}{\max\{1,T/N \cdot \left| a(Y_n^N(\w)) +
b(Y_n^N(\w))\left(W_{\frac{(n+1)T}{N}}(\w) - W_{\frac{nT}{N}}(\w)\right)\right|\}},
\eeqq
for every $n\in\{0,1,...,N-1\}, N\in\bbN$ and all $\w\in\W.$ Η (\ref{eq1000}) is explicit, does not explode and converges strongly to the exact solution $x_t$ of SDE (\ref{eq1}), i.e.,
\beqq\label{eq1013}
\lim_{\Nto}\left(\sup_{0\leq t\leq T}\bfE\Big|x_t-\ov{Y}^N_t\Big|^q\right)=0,
\eeqq
for some $q>0,$ where $\ov{Y}^N_t:=(n+1-\frac{tN}{T})Y_n^N + (\frac{tN}{T} - n)Y_{n+1}^N$ are continuous versions of (\ref{eq1000}) through linear interpolation. It still does not preserve positivity.

For the aforementioned reasons there is an interest in the
construction of numerical schemes to simulate the corresponding
SDEs, that have the desired properties. An attempt to this direction has been made by the first author in (\cite{halidias:2012}, \cite{halidias:2013})
 suggesting the Semi-Discrete method (where, briefly saying, we discretize a part of the SDE). Using this method in  (\cite{halidias:2012}) the author produced a new numerical scheme (but not unique in this situation) for the first aforementioned
problem and proves the strong convergence of the scheme in mean square sense. Later on, in (\cite{halidias:2013}), the author generalizes the idea
of the Semi-Discrete method and uses this generalization to
approximate a class  of  super linear problems, suggesting a
new numerical scheme that preserves positivity in that case, proving again the strong
convergence in the mean square sense.

A basic feature of the Semi-Discrete method is that it is
explicit, compared to other interesting, but implicit methods
(\cite{mao_szpruch:2013b},\cite{mao_szpruch:2013a}), and converges
strongly in the mean square sense to the exact solution of the
original SDE. Moreover, the Semi-Discrete method preserves
positivity  (\cite[Section 3]{halidias:2012}) and it does not
explode in some super-linear problems (\cite[Section
3]{halidias:2013}). The purpose of this paper is to generalize further the method to include non-autonomous coefficients, $a(t,x), b(t,x)$ in (\ref{eq1}) and cover cases like that  of the Heston $3/2$-model.
\section{The setting and the main result.}\label{sec:s1,5}
\setcounter{equation}{0}
\fbox{\begin{minipage}[t] {16cm}
\textbf{Assumption A}  Let $f(s,r,x,y), g(s,r,x,y):[0,T]^2\times \bbR^2\mapsto\bbR$ be such that $f(s,s,x,x)=a(s,x), g(s,s,x,x)=b(s,x),$ where $f,g$ satisfy the following conditions
\beao
|f(s_1,r_1,x_1,y_1) - f(s_2,r_2,x_2,y_2)|&\leq& C_R \left(|s_1-s_2| + |r_1-r_2| + |x_1-x_2| + |y_1-y_2|\right)\\
|g(s_1,r_1,x_1,y_1) - g(s_2,r_2,x_2,y_2)|&\leq& C_R \left(|s_1-s_2| + |r_1-r_2|+|x_1-x_2| + |y_1-y_2|\right.\\
&& + \left.\sqrt{|x_1-x_2|}\right),
\eeao
for  any $R>0$ such that $|x_1|\vee|x_2|\vee|y_1|\vee|y_2|\leq R,$ where the constant $C_R$ depends on $R$ and $x\vee y$ denotes the maximum of $x, y.$\end{minipage}}

Let the equidistant partition $0=t_0<t_1<...<t_N=T$ and $\D=T/N.$
We propose the following Semi-Discrete numerical scheme
\beqq\label{eq1.1} y_t=y_n + \int_{t_n}^{t} f(t_n,s,y_{t_n},y_s)ds
+ \int_{t_n}^{t} g(t_n,s,y_{t_n},y_s)dW_s,\quad t\in[t_n,
t_{n+1}], \eeqq 
where we assume that for every $n\leq N-1,$ (\ref{eq1.1}) has a unique strong solution and $y_n=y_{t_n}, y_0=x_0,$ a.s. In order to compare with the exact solution $x_t,$
which is a continuous time process, we consider the following
interpolation process of the Semi-Discrete approximation, in a
compact form, 
\beqq\label{eq2} y_t=y_0 + \int_{0}^{t}f(\hat{s},s,y_{\hat{s}},y_s)ds +
\int_{0}^{t}g(\hat{s},s,y_{\hat{s}},y_s) dW_s, \eeqq 
where $\hat{s}=t_{n},$ when $s\in[t_{n},t_{n+1}).$ The first and third
variable in $f,g$ denote the discretized part of the original SDE.
We observe from (\ref{eq2}) that in order to solve for $y_t$, we
have to solve an SDE and not an algebraic equation, thus in this
context, we cannot reproduce implicit schemes, but we can
reproduce the Euler scheme if we choose $f(s,r,x,y)=a(s,x)$ and
$g(s,r,x,y)=b(s,x).$

The numerical scheme (\ref{eq2}) converges to the true solution $x_t$ of SDE (\ref{eq1}) and this is stated in the following, which  is our main result.
\bth \label{t1}
Suppose Assumption A holds and (\ref{eq1.1}) has  a unique strong solution for every $n\leq N-1,$ where $x_0\in \bbl^p(\Omega,\bbR), x_0>0$ a.s. Let also
$$
\bfE(\sup_{0\leq t\leq T}|x_t|^p) \vee \bfE(\sup_{0\leq t\leq T}|y_t|^p)<A,
$$
for some $p>2$ and $A>0.$ Then the Semi-Discrete numerical scheme (\ref{eq2}) converges to the true solution of (\ref{eq1}) in the mean square sense, that is
\beqq \label{eq600}
\lim_{\D\rightarrow0}\bfE\sup_{0\leq t\leq T}|y_t-x_t|^2=0.
\eeqq
\ethe

In (\cite{halidias:2013}) the case with no square root term is
treated, thus Theorem \ref{t1} is a generalization of
(\cite[Theorem 1]{halidias:2013}). Section \ref{sec:s2} provides all
the necessary and finally the proof of Theorem \ref{t1}.
Section \ref{sec:s3} gives applications to super linear drift and
diffusion problems with non negative solution, one of which includes the Heston $3/2$-model.  
Section \ref{sec:s5} shows experimentally the order of convergence of the SD method applied to the Heston $3/2$-model. 
The Semi-Discrete scheme is strongly convergent in the mean square sense and preserves positivity of the solution.

\section{Proof of Theorem \ref{t1}.}\label{sec:s2}
\setcounter{equation}{0}
We denote the indicator function of a set $A$ by $\bbi_{A}.$ The constant $C_R$ may vary from line to line and it may depend apart from $R$ on other quantities, like time $T$ for example, which are all constant, in the sense that we don't let them grow to infinity.

\subsection{Error bound for the explicit Semi-Discrete scheme}

\ble \label{l1} Let the assumption of Theorem \ref{t1} hold. Let
$R>0,$ and set the stopping time $\theta_R=\inf\{t\in [0,T]:
|y_t|>R \, \hbox { or } \,  |y_{\hat{t}}|>R\}.$ Then the following
estimate holds \beqq\label{eq3}
\bfE|y_{s\wedge\theta_R}-y_{\wh{s\wedge\theta_R}}|^2 \leq C_{R}
\D, \eeqq where  $C_R$ does not depend on $\D,$ implying
$\sup_{s\in[t_{n_s},t_{n_s+1}]}\bfE|y_{s\wedge\theta_R}-y_{\wh{s\wedge\theta_R}}|^2=O(\D),$
as $\D\downarrow0.$ \ele

\bpf[Proof of Lemma \ref{l1}]
Let $n_s$ integer such that $s\in[t_{n_s},t_{n_s+1}).$ It holds that
\beao
&&|y_{s\wedge\theta_R}-y_{\wh{s\wedge\theta_R}}|^2=\left| \int_{t_{\wh{n_s\wedge\theta_R}}}^{s\wedge\theta_R}f(\hat{u},u,y_{\hat{u}},y_u)du + \int_{t_{\wh{n_s\wedge\theta_R}}}^{s\wedge\theta_R}g(\hat{u},u,y_{\hat{u}},y_u) dW_u\right|^2\\
&\leq&2\left(\int_{t_{\wh{n_s\wedge\theta_R}}}^{s\wedge\theta_R}f(\hat{u},u,y_{\hat{u}},y_u)du\right)^2 + 2\left(\int_{t_{\wh{n_s\wedge\theta_R}}}^{s\wedge\theta_R}g(\hat{u},u,y_{\hat{u}},y_u) dW_u\right)^2\\
&\leq&2\D\int_{t_{\wh{n_s\wedge\theta_R}}}^{s\wedge\theta_R}f^2(\hat{u},u,y_{\hat{u}},y_u)du + 2\left(\int_{t_{\wh{n_s\wedge\theta_R}}}^{s\wedge\theta_R}g(\hat{u},u,y_{\hat{u}},y_u) dW_u\right)^2\\
&\leq&C_R\D^2 +
2\left(\int_{t_{\wh{n_s\wedge\theta_R}}}^{s\wedge\theta_R}g(\hat{u},u,y_{\hat{u}},y_u)
dW_u\right)^2, \eeao where we have used Cauchy-Schwarz inequality
and Assumption A for the function $f.$\footnote{By the fact that we want the problem (\ref{eq1}) to be well posed and by the conditions on $f$ and $g$ we get that $f,g$ are bounded on bounded intervals.} Taking expectations  in the above inequality gives
\beao
\bfE|y_{s\wedge\theta_R}-y_{\wh{s\wedge\theta_R}}|^2 &\leq& C_R\D^2  + 8\bfE\int_{t_{\wh{n_s\wedge\theta_R}}}^{ t_{n_s+1}\wedge\theta_R }g^2(\hat{u},u,y_{\hat{u}},y_u) du\\
&\leq& C_R\D^2  +  C_R \D, \eeao
where in the first step we have used Doob's martingale inequality (\cite[Theorem
1.3.8]{karatzas_shreve:1988}) on the diffusion term,  in the second
step Assumption A for the function $g.$ 
Thus,
$$
\lim_{\D\downarrow0}\frac{\sup_{s\in[t_{n_s},t_{n_s+1}]}\bfE|y_{s\wedge\theta_R}-y_{\wh{s\wedge\theta_R}}|^2}{\D}\leq
C_R,
$$
which justifies the $O(\D)$ notation, (see for example \cite{olver:1997}).
\epf

\subsection{Convergence of the Semi-Discrete scheme in $\bbl^1$}
\bpr\label{pr1}
Let the assumptions of Theorem \ref{t1} hold. Let $R>0,$ and set the stopping time $\theta_R=\inf\{t\in [0,T]: |y_t|>R \, \hbox { or } \,  |x_{t}|>R\}.$ Then we have
\beqq \label{eq4.4}
\sup_{0\leq t\leq T}\bfE|y_{t\wedge\theta_R}-x_{t\wedge\theta_R}|\leq \left[ \left(C_R + \frac{C_R}{me_m}\right)\sqrt{\D} + \left(\frac{C_R}{me_m}+C_R\right)\D + \frac{C_R}{me_m}\D^2 + \frac{C_R}{m} +  e_{m-1}\right]e^{a_{R,m}T},
\eeqq
for any $m>1,$ where
$$
e_m=e^{-m(m+1)/2}, \quad a_{R,m}:=C_R + \frac{C_R}{m}
$$
and $C_R$ does not depend on $\D.$ It holds that $\lim_{m\uparrow\infty}e_m=0.$
\epr

\bpf[Proof of Proposition \ref{pr1}]
Let the non increasing sequence $\{e_m\}_{m\in\bbN}$ with $e_m=e^{-m(m+1)/2}$ and $e_0=1.$ We introduce the following sequence of smooth approximations of $|x|,$ (method of Yamada and Watanabe, \cite{yamada_watanabe:1971})
$$
\phi_m(x)=\int_0^{|x|}dy\int_0^{y}\psi_m(u)du,
$$
where the existence of the continuous function $\psi_m(u)$ with $0\leq \psi_m(u) \leq 2/(mu)$ and support in $(e_m,e_{m-1})$ is justified by $\int_{e_m}^{e_{m-1}}(du/u)=m.$ The following relations hold for $\phi_m\in\bbc^2(\bbR,\bbR)$ with $\phi_m(0)=0,$
 $$
 |x| - e_{m-1}\leq\phi_m(x)\leq |x|, \quad |\phi_{m}^{\prime}(x)|\leq1, \quad x\in\bbR, $$
 $$
 |\phi_{m}^{\prime \prime }(x)|\leq\frac{2}{m|x|}, \,\hbox{ when }  \,e_m<|x|<e_{m-1} \,\hbox{ and }  \,  |\phi_{m}^{\prime \prime }(x)|=0 \,\hbox{ otherwise. }
 $$
We have that
\beqq\label{eq5}
\bfE|y_{t\wedge\theta_R}-x_{t\wedge\theta_R}| \leq e_{m-1} + \bfE\phi_m(y_{t\wedge\theta_R}-x_{t\wedge\theta_R}).
\eeqq
Applying Ito's formula to the sequence $\{\phi_m\}_{m\in\bbN},$ we get
\beao
&&\phi_m(y_{t\wedge\theta_R}-x_{t\wedge\theta_R})= \int_{0}^{t\wedge\theta_R} \phi_m^{\prime}(y_s-x_s) (f(\hat{s},s,y_{\hat{s}},y_s)-f(s,s,x_s,x_s))ds +M_t\\
&&+ \frac{1}{2}\int_{0}^{t\wedge\theta_R} \phi_m^{\prime\prime}(y_s-x_s) (g(\hat{s},s,y_{\hat{s}},y_s)-g(s,s,x_s,x_s))^2 ds \\
\nonumber&\leq& \int_{0}^{t\wedge\theta_R} C_R\left(|y_{\hat{s}}-x_s | + |y_s-x_s| +|\hat{s}-s|\right) ds + M_t\\
&& + \frac{1}{2}\int_{0}^{t\wedge\theta_R} \frac{2}{m|y_s-x_s|} C_R\left(|y_{\hat{s}}-x_s |^{2} + |y_s-x_s|^2 + |y_{\hat{s}}-x_s| + |\hat{s}-s|^2\right) ds \\
&\leq& C_R\int_{0}^{t\wedge\theta_R}|y_s-y_{\hat{s}}|ds + C_R \int_{0}^{t\wedge\theta_R} |y_s-x_s|ds +  C_R \int_{0}^{t\wedge\theta_R} |\hat{s}-s|ds + M_t\\
&&+  \frac{C_R}{m}\int_{0}^{t\wedge\theta_R}\frac{2|y_s - y_{\hat{s}}|^2 + 3|y_s-x_s|^2 + |y_{\hat{s}}-x_s | + |\hat{s}-s|^2}{|y_s-x_s|}ds\\
&\leq& (C_R + \frac{C_R}{me_m})\int_{0}^{t\wedge\theta_R}|y_s-y_{\hat{s}} |ds +
\frac{C_R}{me_m}\int_{0}^{t\wedge\theta_R}|y_s-y_{\hat{s}}|^2ds + (C_R + \frac{C_R}{m}) \int_{0}^{t\wedge\theta_R} |y_s-x_s|ds\\
&&+  \frac{C_R}{me_m}\sum_{k=0}^{[t/\D-1]}\int_{t_k}^{t_{k+1}\wedge\theta_R}|t_k-s|^2ds + C_R\sum_{k=0}^{[t/\D-1]}\int_{t_k}^{t_{k+1}\wedge\theta_R}|t_k-s|ds + \frac{C_R}{m} + M_t\\
&\leq& (C_R + \frac{C_R}{me_m})\int_{0}^{t\wedge\theta_R}|y_s-y_{\hat{s}} |ds +
\frac{C_R}{me_m}\int_{0}^{t\wedge\theta_R}|y_s-y_{\hat{s}}|^2ds\\
&&+ (C_R + \frac{C_R}{m}) \int_{0}^{t\wedge\theta_R} |y_s-x_s|ds + \frac{C_R}{m} + \frac{C_R}{me_m}\D^2 + C_R\D + M_t,
\eeao
where in the second step we have used Assumption A for the functions $f, g$ and the properties of $\phi_m$ and
$$
M_t:= \int_{0}^{t\wedge\theta_R} \phi_m^{\prime}(y_u-x_u)
(g(\hat{u},u,y_{\hat{u}},y_u)-g(u,u,x_u,x_u)) dW_u.
$$
Taking expectations in the above inequality yields
\beao
\bfE\phi_m(y_{t\wedge\theta_R}-x_{t\wedge\theta_R})&\leq& \left(C_R + \frac{C_R}{me_m}\right) \int_{0}^{t\wedge\theta_R}\bfE|y_s-y_{\hat{s}}|ds + \left(C_R + \frac{C_R}{m}\right) \int_{0}^{t\wedge\theta_R} \bfE|y_s-x_s|ds\\
&&+ \frac{C_R}{me_m}\int_{0}^{t\wedge\theta_R}\bfE|y_s-y_{\hat{s}}|^2ds + \frac{C_R}{m} +   \frac{C_R}{me_m}\D^2 + C_R\D + \bfE M_t\\
&\leq& \left(C_R + \frac{C_R}{me_m}\right) \sqrt{\D} + \left(\frac{C_R}{me_m}+C_R\right)\D  + \frac{C_R}{me_m}\D^2 + \frac{C_R}{m}\\
&& + \left(C_R + \frac{C_R}{m}\right) \int_{0}^{t\wedge\theta_R} \bfE|y_s-x_s|ds,
\eeao
where we have used  Lemma \ref{l1} and the fact that $\bfE M_t=0$.\footnote{The function $h(u)=\phi_m^{\prime}(y_u-x_u)
(g(\hat{u},u,y_{\hat{u}},y_u)-g(u,u,x_u,x_u))$ belongs to the space $\bbm^2([0,t\wedge\theta_R];\bbR)$ of real valued measurable $\bbf_t-$adapted processes such that $\bfE\int_0^{t\wedge\theta_R}|h(u)|^2du<\infty$ thus (\cite[Theorem 1.5.8]{mao:1997}) implies $\bfE M_t=0$.} Thus (\ref{eq5}) becomes
\beao
&&\bfE|y_{t\wedge\theta_R}-x_{t\wedge\theta_R}| \leq \left(C_R + \frac{C_R}{me_m}\right) \sqrt{\D} + \left(\frac{C_R}{me_m}+C_R\right)\D + \frac{C_R}{me_m}\D^2 +\frac{C_R}{m} +  e_{m-1}\\
&& + \left(C_R + \frac{C_R}{m}\right) \int_{0}^{t\wedge\theta_R} \bfE|y_s-x_s|ds\\
&\leq&\left[\left(C_R + \frac{C_R}{me_m}\right) \sqrt{\D} +
\left(\frac{C_R}{me_m}+C_R\right)\D  +\frac{C_R}{me_m}\D^2 + \frac{C_R}{m} +  e_{m-1}\right]e^{a_{R,m}t},
\eeao
where in the last step we have used the Gronwall inequality (\cite[Relation 7]{gronwall:1919}) and $a_{R,m}=C_R + \frac{C_R}{m}.$ Taking the supremum over all $0\leq t \leq T$
gives (\ref{eq4.4}).
\epf

\subsection{Convergence of the Semi-Discrete scheme in $\bbl^2$}

Set the stopping time $\theta_R=\inf\{t\in [0,T]: |y_t|>R \, \mbox{or} \, |x_t|>R\},$ for some $R>0$ big enough. We have that
\beam
\nonumber
\bfE\sup_{0\leq t\leq T}|y_t-x_t|^2\!\!\!&=&\!\!\!\bfE\sup_{0\leq t\leq T}|y_t-x_t|^2\bbi_{(\theta_R> t)} +
\bfE\sup_{0\leq t\leq T}|y_t-x_t|^2\bbi_{(\theta_R\leq t)}\\
\nonumber\!\!\!&\leq&\!\!\!\bfE\sup_{0\leq t\leq T }|y_{t\wedge\theta_R}-x_{t\wedge\theta_R}|^2 +
\frac{2\del}{p}\bfE\sup_{0\leq t\leq T}|y_t - x_t|^{p} + \frac{(p-2)}{p\del^{2/(p-2)}}\bfP(\theta_R\leq T)\\
\nonumber\!\!\!&\leq&\!\!\!\bfE\sup_{0\leq t\leq T }|y_{t\wedge\theta_R}-x_{t\wedge\theta_R}|^2 + \frac{2^p\del}{p}\bfE\sup_{0\leq t\leq T}\left(|y_t|^p + |x_t|^{p}\right) + \frac{(p-2)}{p\del^{2/(p-2)}}\bfP(\theta_R\leq T)\\
\label{eq8}\!\!\!&\leq&\!\!\!\bfE\sup_{0\leq t\leq T }|y_{t\wedge\theta_R}-x_{t\wedge\theta_R}|^2 + \frac{2^{p+1}\del A}{p} + \frac{(p-2)}{p\del^{2/(p-2)}}\bfP(\theta_R\leq T),
\eeam
where in the second step we have applied Young inequality,
$$
ab \leq \frac{\del}{r}a^r + \frac{1}{q\del^{q/r}}b^q,
$$
for $a=\sup_{0\leq t\leq T}|y_t - x_t|^2, b=\bbi_{(\theta_R\leq t)}, r=p/2, q=p/(p-2)$ and $\del>0,$ in the third step we have used the elementary inequality $(\sum_{i=1}^n a_i)^p \leq n^{p-1}\sum_{i=1}^n a_i^p,$ with $n=2,$ and $A$ comes from the moment bound assumption. It holds that
$$
\!\bfP(\theta_R\leq T)\!\leq\!\bfE\!\left(\bbi_{(\theta_R\leq
T)}\frac{|y_{\theta_R}|^{p}}{R^p}\right) + 
\bfE\!\left(\bbi_{(\theta_R\leq T)}\frac{|x_{\theta_R}|^{p}}{R^p}\right)\!\leq\!
\frac{1}{R^p}\!\left(\bfE\sup_{0\leq t\leq T}|x_t|^{p} +
\bfE\sup_{0\leq t\leq T}|y_t|^{p}\right)\!\leq\!\frac{2A}{R^p},
$$
thus (\ref{eq8}) becomes
\beqq\label{eq10}
\bfE\sup_{0\leq t\leq T}|y_t-x_t|^2 \leq \bfE\sup_{0\leq t\leq T }|y_{t\wedge\theta_R}-x_{t\wedge\theta_R}|^2 + \frac{2^{p+1}\del A}{p} + \frac{2(p-2)A}{p\del^{2/(p-2)}R^p}.
\eeqq
We estimate the difference $|e_{t\wedge\theta_R}|^2:=|y_{t\wedge\theta_R}-x_{t\wedge\theta_R}|^2.$ It holds that
\beao
&&\!\!\!\!|e_{t\wedge\theta_R}|^2\!\!=\!\!\left|\int_{0}^{t\wedge\theta_R}\!\!\left(f(\hat{s},s,y_{\hat{s}},y_s)-f(s,s,x_s,x_s)\right)ds + \!\! \int_{0}^{t\wedge\theta_R}\!\!\left( g(\hat{s},s,y_{\hat{s}},y_s)-g(s,s,x_s,x_s) \right)\!dW_s\right|^2\\
&\leq& 2T\int_{0}^{t\wedge\theta_R}C_R \left(|y_{\hat{s}}-x_s|^2 + |y_s-x_s|^2+ |\hat{s}-s|^2\right) ds + 2 |M_t|^2\\
&\leq& C_R\int_{0}^{t\wedge\theta_R}|y_s-y_{\hat{s}}|^2ds +  C_R\int_{0}^{t\wedge\theta_R}|y_s-x_s|^2ds + C_R\int_{0}^{t\wedge\theta_R}|\hat{s}-s|^2ds + 2|M_t|^2\\
&\leq& C_R\int_{0}^{t\wedge\theta_R}|y_s-y_{\hat{s}}|^2ds +  C_R\int_{0}^{t\wedge\theta_R}|y_s-x_s|^2ds + C_R\sum_{k=0}^{[t/\D-1]}\int_{t_k}^{t_{k+1}\wedge\theta_R}|t_k-s|^2ds + 2|M_t|^2,
\eeao
where in the second step we have used Cauchy-Schwarz inequality and Assumption A for $f$ and
$$
M_t:=\int_{0}^{t\wedge\theta_R}\left( g(\hat{s},s,y_{\hat{s}},y_s)-g(s,s,x_s,x_s) \right)dW_s.
$$
Taking the supremum over all $t\in[0,T]$ and then expectations we have
\beam
\nonumber
&&\bfE\sup_{0\leq t\leq T}|y_{t\wedge\theta_R}-x_{t\wedge\theta_R}|^2
\leq C_R \bfE\left(\int_{0}^{T\wedge\theta_R}|y_s-y_{\hat{s}}|^2ds\right)+ 2\bfE\sup_{0\leq t\leq T} |M_t|^2\\
\nonumber&& + C_R\int_{0}^{T}\bfE\sup_{0\leq l\leq s}|y_{l\wedge\theta_R}-x_{l\wedge\theta_R}|^2ds + C_R\D^2\\
\label{eq9}&\leq& C_R\int_{0}^{T\wedge\theta_R}\bfE|y_s-y_{\hat{s}}|^2ds  + 8\bfE |M_T|^2 + C_R \int_{0}^{T}\bfE\sup_{0\leq l\leq s}|y_{l\wedge\theta_R}-x_{l\wedge\theta_R}|^2ds + C_R\D^2,
\eeam
where in the last step we have used Holder's inequality and Doob's martingale inequality with $p=2,$ since  $M_t$ is an $\bbR-$valued martingale that belongs to $\bbl^2.$ It holds that
\beao
&&\bfE |M_T|^2:=\bfE\left|\int_{0}^{T\wedge\theta_R} \left( g(\hat{s},s,y_{\hat{s}},y_s)-g(s,s,x_s,x_s) \right) dW_s\right|^2\\
&=& \bfE\left(\int_{0}^{T\wedge\theta_R} \left( g(\hat{s},s,y_{\hat{s}},y_s)-g(s,s,x_s,x_s) \right) ^2ds\right)\\
&\leq& C_R\bfE\left(\int_{0}^{T\wedge\theta_R} \left( |y_{\hat{s}}- x_s|^{2} + |y_s-x_s|^2 + |y_{\hat{s}}- x_s| + |\hat{s}-s|^2 \right) ds\right)\\
&\leq& C_R\int_{0}^{T\wedge\theta_R} \bfE|y_s-y_{\hat{s}}|^2ds + C_R\int_{0}^{T}\bfE\sup_{0\leq l\leq s}|y_{l\wedge\theta_R}-x_{l\wedge\theta_R}|^2ds + C_R\int_{0}^{T\wedge\theta_R} \bfE|y_{\hat{s}}-x_s|ds + C_R\D^2,
\eeao
where we have used Assumption A for $g.$ Relation (\ref{eq9}) becomes
 \beao
&& \bfE\sup_{0\leq t\leq T}|y_{t\wedge\theta_R}-x_{t\wedge\theta_R}|^2\leq C_R\int_0^{T\wedge\theta_R} \bfE|y_s-y_{\hat{s}}|^2ds  + C_R\int_{0}^{T}\bfE\sup_{0\leq l\leq s}|y_{l\wedge\theta_R}-x_{l\wedge\theta_R}|^2ds\\
&&+ C_R\int_{0}^{T\wedge\theta_R} \left(\bfE|y_s - y_{\hat{s}}| + \bfE|y_s - x_s|\right)ds + C_R\D^2\\
 &\leq& C_{R}\sqrt{\D} + C_{R}\D + C_R\D^2 + C_R\int_{0}^{T}\bfE\sup_{0\leq l\leq s}|y_{l\wedge\theta_R}-x_{l\wedge\theta_R}|^2ds + C_R\int_{0}^{T\wedge\theta_R} \bfE|y_s - x_s|ds,
 \eeao
where we have used Lemma \ref{l1} and Jensen's inequality for the concave function $\phi(x)=\sqrt{x}.$  The integrand of the last term is bounded, from Proposition \ref{pr1}, by
$$
K_{R,\D,m}(s)\!:=\!\left[ \left(C_R + \frac{C_R}{me_m}\right)\!\sqrt{\D} + \left(\frac{C_R}{me_m}+C_R\right)\D + \frac{C_R}{me_m}\D^2 + \frac{C_R}{m} + e_{m-1}\right] e^{a_{R,m}s},
$$
where $s\in[0,T\wedge \theta_R].$ Application of the Gronwall inequality implies
 $$
\bfE\sup_{0\leq t\leq T}|y_{t\wedge\theta_R}-x_{t\wedge\theta_R}|^2 \leq \left(C_R\sqrt{\D} + C_R\D + C_R K_{R,\D,m}(T)\right)e^{C_R}\leq C_{R,\D,m}.
$$
Note that, given $R>0,$ the quantity $C_{R,\D,m}$ can be arbitrarily small by choosing big enough $m$ and small enough  $\D.$ Relation (\ref{eq10}) becomes,
 \beao
 \bfE\sup_{0\leq t\leq T}|y_t-x_t|^2 &\leq& C_{R,\D,m} +  \frac{2^{p+1}\del A}{p} + \frac{2(p-2)A}{p\del^{2/(p-2)}R^p}\\
  &&:=I_1 + I_2 + I_3.
\eeao

Given any $\ep>0,$ we may first choose $\del$ such that
$I_2<\ep/3,$ then choose $R$ such that  $I_3<\ep/3,$ then  $m>1$
and finally $\D$ such that $I_1<\ep/3$ concluding $\bfE\sup_{0\leq
t\leq T}|y_t-x_t|^2 <\ep$ as required to verify (\ref{eq600}).

\section{Superlinear examples.}\label{sec:s3}
\setcounter{equation}{0}

\subsection{Example I}\label{ssec:3.1}

We study the numerical approximation of the following SDE,
\beqq\label{eq11} x_t=x_0 + \int_{0}^{t}(k_1(s) x_s -k_2
(s)x_s^{2})ds + \int_{0}^{t} k_3(s)x_s^{3/2} \phi(x_s) dW_s,\quad
t\in [0,T], \eeqq where $\phi(\cdot)$ is a locally Lipschitz and
bounded function with locally Lipschitz constant $C_R^\phi,$
bounding constant $K_\phi$, $x_0$ is independent of all
$\{W_t\}_{0\leq t\leq T}, x_0\in \bbl^{4p}(\Omega,\bbR)$ for some
$2<p$ and $x_0>0,$ a.s., $\mathbb{E}(x_0)^{-2} < A$, $k_1(\cdot),
k_2(\cdot), k_3(\cdot)$ are positive and bounded functions with
$k_{2,\min}> \frac{7}{2}(K_\phi k_{3,\max})^2.$  Model
(\ref{eq11}) has super linear drift and diffusion coefficients.

We propose the following Semi-Discrete numerical scheme
\beqq\label{eq12}
y_t=y_n + \int_{t_n}^{t} (k_1(s) - k_2(s)y_{t_n})y_sds   + \int_{t_n}^{t} k_3(s)\sqrt{y_{t_n}}\phi(y_{t_n})y_s  dW_s,\quad t\in[t_n, t_{n+1}],
\eeqq
where $y_n=y_n(t_n),$  for $n\leq T/\D$ and $y_0=x_0,$ a.s., or in a  more compact form,
\beqq\label{eq13}
y_t=y_0 + \int_{0}^{t}(k_1(s) - k_2(s)y_{\hat{s}})y_sds +
\int_{0}^{t}k_3(s)\sqrt{y_{\hat{s}}}\phi(y_{\hat{s}})y_s dW_s,
\eeqq
where $\hat{s}=t_n,$ when $s\in[t_n,t_{n+1}).$ The linear SDE (\ref{eq13}) has a solution which, by use of Ito's formula, has the explicit form
\beqq\label{eq14}
y_t=x_0\exp\Big\{\int_{0}^{t}\left(k_1(s) - k_2(s) y_{\hat{s}} - k_3^2(s)\frac{y_{\hat{s}}\phi^2(y_{\hat{s}})}{2}\right) ds +
\int_{0}^{t}k_3(s)\sqrt{y_{\hat{s}}}\phi (y_{\hat{s}}) dW_s\Big\},
\eeqq
where $y_t=y_t(t_0,x_0).$
\bpr \label{pr0}
The Semi-Discrete   numerical scheme (\ref{eq13}) converges to the true solution of
(\ref{eq11}) in the mean square sense, that is
\beqq \label{eq15}
\lim_{\D\rightarrow0}\bfE\sup_{0\leq t\leq T}|y_t-x_t|^2=0.
\eeqq
\epr

\subsubsection{Proof of Proposition \ref{pr0}}\label{ssec:s40}

In order to prove Proposition \ref{pr0} we need to verify the assumptions of Theorem \ref{t1}. Let
\beao
a(s,x)=k_1(s) x -k_2(s) x^{2}, && f(s,r,x,y)=(k_1(s)-k_2(s)x)y,\\
b(s,x)= k_3(s)x^{3/2}\phi(x),  && g(s,r,x,y)=k_3(s)\sqrt{x}\phi(x)y.
\eeao
We verify Assumption A for $f.$ Let $R>0$ such that $|x_1|\vee|x_2|\vee|y_1|\vee|y_2|\vee|s|\vee|r|\leq R.$ We have that
\beao
&&|f(s,r,x_1,y_1)-f(s,r,x_2,y_2)|=|(k_1(s)-k_2(s) x_1)y_1 - (k_1(s)-k_2(s)x_2)y_2|\\
&\leq& |k_{1}(s)||y_1-y_2| + |k_{2}(s)|(|x_2||y_1-y_2| + |y_1||x_1-x_2|)\\
&\leq& (|k_{1,\max}| + |k_{2,\max}|R)|y_1-y_2| + |k_{2,\max}|R|x_1-x_2|\\
&\leq& C_R\left( |x_1-x_2| + |y_1-y_2| \right),
\eeao
thus, Assumption A holds for $f$ with $C_R:=|k_{1,\max}| + |k_{2,\max}|R.$

We verify Assumption A for $g.$ Let $R>0$ such that $|x_1|\vee|x_2|\vee|y_1|\vee|y_2|\vee|s|\vee|r|\leq R.$ We have that
\beao
&&|g(s,r,x_1,y_1)-g(s,r,x_2,y_2)|=|k_3(s)\sqrt{x_1}\phi (x_1)y_1 - k_3(s)\sqrt{x_2}\phi(x_2)y_2|\\
&\leq& |k_{3}(s)| \left(\sqrt{x_1}| \phi (x_1)| |y_1-y_2| + |y_2| \big| \sqrt{x_1}\phi(x_1) - \sqrt{ x_1}\phi( x_2) + \sqrt{x_1}\phi(x_2) - \sqrt{x_2}\phi(x_2)\big|\right)\\
&\leq& |k_{3,\max}|\left(K_\phi \sqrt{R} |y_1-y_2| + R \sqrt{x_1} | \phi(x_1) - \phi(x_2)| +  R K_\phi|\sqrt{x_1} - \sqrt{x_2}|\right)\\
&\leq& |k_{3,\max}|\left( K_\phi \sqrt{R} |y_1-y_2| + R^{3/2} C^\phi_R| x_1 - x_2| +  R K_\phi \sqrt{|x_1 - x_2|}\right)\\
&\leq& C_R \left(|x_1 -x_2| +|y_1-y_2| + \sqrt{|x_1 -x_2|}\right),
\eeao where we have used the fact that the function $\sqrt{x}$ is
$1/2-$Holder continuous and $C_R:=|k_{3,\max}|\left(C_R^\phi R^{3/2}\vee
K_\phi \sqrt{R}\vee K_\phi R\right).$ Thus, Assumption A holds for
$g$.

\subsubsection{Moment bound for original SDE}

\ble\label{l4} In the previous setting it holds that $x_t > 0$
a.s. \ele

\bpf[Proof of Lemma \ref{l4}] Set the stopping time $\theta_R=\inf\{t\in [0,T]: x_t^{-1}>R\},$ for some $R>0,$ with the convention that $\inf\emptyset=\infty.$ Application of Ito's formula on
$x_{t\wedge\theta_R}^{-2}$ implies,
\beao
&&(x_{t\wedge \theta_R})^{-2} = (x_0)^{-2} + \int_{0}^{t\wedge \theta_R}(-2)(x_s)^{-3}(k_1(s)x_s- k_2(s)x_s^2) ds\\
&&+  \int_{0}^{t\wedge\theta_R}\frac{(-2)(-3)}{2}(x_s)^{-4}k_3^2(s)x_s^{3}\phi^2(x_s)ds + \int_{0}^{t\wedge \theta_R}(-2)k_3(s) (x_s)^{-3} x_{s}^{3/2}\phi(x_s)dW_s\\
&\leq&\!\!\!(x_0)^{-2} + \int_{0}^{t\wedge \theta_R}(-2k_1(s) x_s^{-2}  + 2k_2(s)x_s^{-1} + 3k_3^2(s)K^2_\phi x_s^{-1})ds\\
&&+\int_{0}^{t}(-2)k_3(s) x_s^{-3/2}\phi(x_s)\,\bbi_{(0,t\wedge \theta_R)}(s)dW_s\\
&\leq&\!\!\!(x_0)^{-2} + \int_{0}^{t\wedge \theta_R}\!\Big(-2k_1(s) x_s^{-2}  + (2k_2(s)+ 3k_3^2(s)K^2_\phi)\left( x_s^{-1} \bbi_{(0,1]}(x_s) + x_s^{-1} \bbi_{(1,\infty]}(x_s)\right)\Big)ds\\
&&+ M_t\\
&\leq&\!\!\!(x_0)^{-2} + 2k_{2,\max}T+ 3k_{3,\max}^2K^2_\phi T + \int_{0}^{t}(2k_2(s)
+3k_3^2(s)K^2_\phi)x_s^{-2}\,\bbi_{(0,t\wedge \theta_R)}(s)ds + M_t,
\eeao
where
$$
M_t:=\int_{0}^{t}(-2)k_3(s)x_s^{-3/2}\phi(x_s)\,\bbi_{(0,t\wedge
\theta_R)}(s)dW_s.
$$
Taking expectations in the above inequality and using the fact that $\bfE M_t=0,$\footnote{The function $h(u)=(-2)k_3(u) x_u^{-3/2}\phi(x_u)\,\bbi_{(0,t\wedge
\theta_R)}(u)$ belongs to the space $\bbm^2([0,t];\bbR)$ thus (\cite[Theorem 1.5.8]{mao:1997}) implies $\bfE M_t=0$.} we get that
\beao
\bfE(x_{t\wedge \theta_R}^{-2})\!\!\!\!&\leq&\!\!\!\!\bfE(x_0)^{-2} + 2k_{2,\max}T+ 3(k_{3,\max}K_\phi)^2T + (2k_{2,\max} +3(k_{3,\max}K_\phi)^2)\int_{0}^{t}\bfE (x_{s\wedge \theta_R})^{-2}ds\\
&\leq&\!\!\!\!\left(\bfE(x_0)^{-2} + 2k_{2,\max}T+ 3k_{3,\max}^2K^2_\phi T
\right)e^{(2k_{2,\max} +3k_{3,\max}^2K^2_\phi)T}<C, \eeao
where we have used Gronwall inequality with $C$ independent of $R.$ We
have that
\beqq\label{eq17}
(x_{t\wedge\theta_R})^{-2}=(x_{\theta_R})^{-2}\bbi_{(\theta_R\leq
t)}+ (x_{t})^{-2}\bbi_{(t < \theta_R)}=R^{2}\bbi_{(\theta_R\leq
t)}+ (x_{t})^{-2}\bbi_{(t<\theta_R)}. \eeqq

By relation (\ref{eq17})  we have that,
$$
\bfE\left(\frac{1}{x_{t\wedge\theta_R}^{2}}\right)=R^2\bfP(\theta_R\leq t) +
\bfE\left(\frac{1}{x_t^{2}}\bbi_{(t<\theta_R)}\right)<C,
$$
thus
$$
\bfP(x_t\leq0)=\bfP\left(\bigcap_{R=1}^\infty \Big\{x_t<
\frac{1}{R} \Big\}\right)=\lim_{\Rto}\bfP\left( \Big\{x_t<
\frac{1}{R} \Big\}\right)\leq  \lim_{\Rto}\bfP(\theta_R\leq t)=0.
$$
We conclude that $x_t > 0$ a.s. \epf

\ble\label{l2}
In the previous setting it holds that
$$
\bfE(\sup_{0\leq t\leq T}(x_t)^p)<A_1,
$$
for some $A_1>0$ and any $2< p \leq k_{2,min}/(K_{\phi}k_{3,max})^2.$ \ele 
\bpf[Proof of Lemma \ref{l2}] In the case of $x$'s
outside a finite ball of radius $R,$ with $R>1,$ and $s\in[0,T]$
we have that \beao
J(s,x)&:=&\frac{xa(s,x) + (p-1)b^2(s,x)/2}{1+x^2}=\frac{x(k_1(s) x-k_2(s)x^{2}) + (p-1)k_3^2(s)[x^{3/2}\phi(x)]^2/2}{1+x^2}\\
&=&\frac{k_1(s)x^2 - k_2(s) x^{3} + 0.5(p-1)k_3^2(s)x^{3}\phi^2(x)}{1+x^2}\\
&\leq&\frac{k_{1,\max}x^2 + \Big(0.5(p-1)(k_{3,\max}K_\phi)^2 -
k_{2,\min}\Big)x^{3}}{1+x^2}\leq k_{1,\max}, \eeao where the last
inequality is valid for all $p$ such that $p \leq 1 +
2k_{2,\min}/(K_\phi k_{3,\max})^2.$ Thus $J(s,x)$ is bounded for
all $(s,x)\in[0,T] \times\bbR,$ since when $|x|\leq R$ we have
that $J(s,x)$ is finite and say $J(s,x)\leq C.$ Since $C$ is
positive, application of (\cite[Theorem 2.4.1]{mao:1997}) implies
$$
\bfE (x_t)^p\leq 2^{(p-2)/2}(1+\bfE (x_0)^p)e^{Cpt},
$$
for any $2<p \leq 1 + 2k_{2,\min}/(K_\phi k_{3,\max})^2$ and all $t\in[0,T]$.
Using Ito's formula on $(x_t)^p,$ with $p \leq k_{2,\min}/(K_\phi k_{3,\max})^2$
(in order to use Doob's martingale inequality later) we have that
\beao
&&(x_{t})^p = (x_0)^p + \int_{0}^{t}p(x_s)^{p-1}(k_1(s)x_s- k_2(s)x_s^2) ds\\
&&+  \int_{0}^{t}\frac{p(p-1)}{2}(x_s)^{p-2}[k_3(s)x_s^{3/2}\phi (x_s)]^2 ds + \int_{0}^{t}pk_3(s)(x_s)^{p-1}x_s^{3/2}\phi (x_s)dW_s\\
&\leq& (x_0)^p + p\int_{0}^{t}\left[k_1(s)(x_s)^{p} + \left(\frac{p-1}{2}k_{3,\max}^2K^2_\phi-k_2\right)(x_s)^{p+1}\right]ds + M_t\\
&\leq& (x_0)^p + p\int_{0}^{t}k_1(s)(x_s)^{p}ds + M_t,
\eeao
where $M_t=\int_{0}^{t}pk_3(s)\phi (x_s)(x_s)^{p+1/2}dW_s.$ Taking the supremum and then expectations in the above inequality we get
\beao
\bfE(\sup_{0\leq t\leq T}(x_{t})^p) &\leq& \bfE(x_0)^p + pk_{1,\max}\bfE\left(\sup_{0\leq t\leq T}\int_{0}^{t} (x_s)^{p}ds\right) + \bfE\sup_{0\leq t\leq T}M_t\\
&\leq&\bfE(x_0)^p + pk_{1,\max}\int_{0}^{T} \bfE(\sup_{0\leq l\leq s}(x_l)^{p})ds + \sqrt{\bfE\sup_{0\leq t\leq T}M_t^2}\\
&\leq&\left(\bfE(x_0)^p + \sqrt{4\bfE M_T^2}\right)e^{pk_{1,\max}T}:=A_1,
\eeao
where in the last step we have used Doob's martingale inequality to the diffusion term $M_t$\footnote{The function $h(u)=pk_3(u)\phi (x_u)(x_u)^{p+1/2}$ belongs to the family $\bbm^2([0,T];\bbR)$ thus (\cite[Theorem 1.5.8]{mao:1997}) implies $\bfE M_t^2=\bfE(\int_0^t h(u)dW_u)^2=\bfE \int_0^t h^2(u)du,$ i.e. $M_t\in\bbl^2(\W;\bbR)$.}  and Gronwall inequality.
\epf

\subsubsection{Moment bound for Semi-Discrete approximation}

\ble\label{l3}
In the previous setting it holds that
$$\bfE(\sup_{0\leq t\leq T}(y_t)^p)<A_2,
$$
for some $A_2>0$ and for any $2<p \leq 1/4 + \frac{k_{2,min}}{2
(k_{3,max} K_{\phi})^2}.$ \ele 
\bpf[Proof of Lemma \ref{l3}] Set the stopping time
$\theta_R=\inf\{t\in [0,T]: y_t>R \},$ for some $R>0,$ with the
convention that $\inf\emptyset=\infty.$ Application of Ito's
formula on $(y_{t\wedge \theta_R})^q,$ with $q=4p$ implies, \beao
&&(y_{t\wedge \theta_R})^q = (y_0)^q + \int_{0}^{t\wedge \theta_R}q(y_s)^{q-1}(k_1(s)- k_2(s) y_{\hat{s}})y_s ds\\
&&+  \int_{0}^{t\wedge\theta_R}\frac{q(q-1)}{2}(y_s)^{q-2}\left[ k_3(s)\sqrt{y_{\hat{s}}}\phi (y_{\hat{s}})y_s\right]^2 ds + \int_{0}^{t\wedge \theta_R}qk_3(s)(y_s)^{q-1}\sqrt{y_{\hat{s}}}\phi (y_{\hat{s}})y_sdW_s\\
&=&(x_0)^q + \int_{0}^{t\wedge \theta_R}\left(q(k_1(s)- k_2(s) y_{\hat{s}}) + \frac{q(q-1)k_3^2(s)}{2}y_{\hat{s}}\phi^2 (y_{\hat{s}})\right)(y_s)^{q}ds\\
&& + \int_{0}^{t\wedge \theta_R}q k_3(s)\sqrt{y_{\hat{s}}}\phi (y_{\hat{s}}) (y_s)^{q}dW_s\\
&\leq& (x_0)^q + q\int_{0}^{t}\left[k_1(s) + \left(\frac{q-1}{2}k_{3,\max}^2K^2_\phi-k_{2,\min}\right)y_{\hat{s}}\right](y_s)^{q}\,\bbi_{(0,t\wedge \theta_R)}(s)ds + M_t\\
&\leq& (x_0)^q + q\int_{0}^{t}k_1(s) (y_s)^{q}\,\bbi_{(0,t\wedge
\theta_R)}(s)ds + M_t, \eeao
where the last inequality is valid for $q\leq 1 + 2k_{2,\min}/(k_{3,\max}K_\phi)^2$  and $$M_t:=\int_{0}^{t\wedge \theta_R}q k_3(s)\sqrt{y_{\hat{s}}}\phi(y_{\hat{s}}) (y_s)^{q}dW_s.$$ Taking expectations and using that $\bfE M_t=0$ we get
$$
\bfE(y_{t\wedge \theta_R})^{q} \leq \bfE(x_0)^{q} +
qk_{1,\max}\int_{0}^{t}\bfE (y_{s\wedge \theta_R})^{q}ds,
$$
Application of the Gronwall inequality implies
$$
\bfE(y_{t\wedge \theta_R})^{q} \leq \bfE(x_0)^{q}e^{qk_{1,\max} T}.
$$
We have that
$$
(y_{t\wedge \theta_R})^{q}=(y_{\theta_R})^{q}\bbi_{(\theta_R\leq
t)}+ (y_{t})^{q}\bbi_{(t < \theta_R)}=R^{q}\bbi_{(\theta_R\leq
t)}+ (y_{t})^{q}\bbi_{(t<\theta_R)},
$$
thus taking expectations in the above inequality and using the
estimated upper bound for $\bfE(y_{t\wedge \theta_R})^{q}$ we
arrive at
$$
\bfE (y_{t})^{q}\bbi_{(t<\theta_R)}\leq \bfE(x_0)^{q}e^{qk_{1,\max} T}
$$
and taking limits in both sides as $\Rto$ we get that
$$
\lim_{\Rto}\bfE (y_{t})^{q}\bbi_{(t<\theta_R)}\leq
\bfE(x_0)^{q}e^{qk_{1,\max} T}.
$$
Fix $t$. The sequence $(y_{t})^{q}\bbi_{(t<\theta_R)}$ is
nondecreasing in $R$ since $\theta_R$ is increasing in $R$ and
$t\wedge \theta_R\rightarrow t$ as $\Rto$ and
$(y_{t})^{q}\bbi_{(t<\theta_R)}\rightarrow(y_{t})^{q}$ as $\Rto,$
thus the monotone convergence theorem
implies
\beqq \label{eq16}
\bfE (y_{t})^{q}\leq
\bfE(x_0)^{q}e^{qk_{1,\max} T}, \eeqq
for any $q \leq 1 + \frac{2k_{2,\min}}{(k_{3,\max}K_\phi)^2}.$ Following the same lines as in Lemma \ref{l2}, i.e. using again Ito's formula on $(y_t)^p$, taking the
supremum and then using Doob's martingale inequality on the
diffusion term we obtain the desired result. Note that in this
last step we need $2 k_{2,\min} > 7(k_{3,\max} K_{\phi})^2$. \epf

\bre\label{r00001}
\begin{itemize}
    \item [(i)] Proposition \ref{pr0} implies that our explicit numerical scheme converges in the mean square sense. Moreover, by (\ref{eq14}) we get that our numerical scheme preserves positivity, which is a desirable modelling property
(\cite{appleby_et_al:2010}, \cite{kahl_et_al:2008}). Example
(\ref{eq11}) covers  the $3/2-$model (\ref{eq1003}), in the case where $\phi(\cdot), k_1(\cdot), k_2(\cdot), k_3(\cdot)$ are constant, and
super-linear problems both in drift and diffusion.

\item [(ii)] Moreover, note that in the analysis that we followed,
we did not discretize the coefficients $k_i$.  In general, by
Theorem \ref{t1},  we are free to discretize any of the
$k_i(\cdot), i=1,2,3,$ functions at any degree. Thus, we can fully
discretize every $k_i(\cdot), i=1,2,3,$ meaning that (\ref{eq12})
will become \beqq\label{eq121} y_t=y_n + \int_{t_n}^{t} (k_1(t_n)
- k_2(t_n)y_{t_n})y_sds   + \int_{t_n}^{t}
k_3(t_n)\sqrt{y_{t_n}}\phi(y_{t_n})y_s  dW_s,\quad t\in[t_n,
t_{n+1}], \eeqq or semi-discretize every  $k_i(\cdot), i=1,2,3,$
\beqq\label{eq122} y_t=y_n + \int_{t_n}^{t} (\hat{k_1}(s,t_n) -
\hat{k_2}(s,t_n)y_{t_n})y_sds   + \int_{t_n}^{t}
\hat{k_3}(s,t_n)\sqrt{y_{t_n}}\phi(y_{t_n})y_s dW_s,\quad
t\in[t_n, t_{n+1}], \eeqq where $\hat{k_i}(t,t)=k_i(t), i=1,2,3.$
The only difference in that situation is  that we require,
$\hat{k}_i(\cdot,\cdot), i=1,2,3$ to be locally Lipschitz in both
variables.

\item [(iii)] One more point of discussion is the dependence on $\w$ that we can assume on the coefficients $k_i$'s. Specifically, we consider the more general SDE
\beqq\label{eq101}
x_t=x_0 + \int_{0}^{t}a_\w(s,x_s)ds + \int_{0}^{t}b_\w(s,x_s)dW_s,\quad t\in [0,T].
\eeqq
Then, assuming that it admits a unique strong solution, our method seems to work. In the example discussed here, an extra condition on the $k_i$'s would be of the form
$$
|k_i(t,\w)|\leq C, t\in[0,T], \w\in\W, i=1,2,3.
$$

\item [(iv)] We illustrate our method in the case $\phi(x) =
\sin(x).$  Then the diffusion term $b(x)$ takes positive and
negative values and thus method (\cite{neuenkirch_szpruch:2012})
does not work since it requires $b(x)>0$ in order to use the
Lamperti-type transformation, as well as Milstein method
(\cite{higham_et_al:2013}) since for the same reason their
Assumption 2.7 is violated. The only method that we know and can
be used for this situation is the Tamed-Euler method
(\cite{hutzenhaler_et_al:2012}, \cite{hutzenhaler_jentzen:2012})
but the drawback is that it does not preserve positivity.

Below, we compare our scheme, in the case where $k_1(\cdot), k_2(\cdot), k_3(\cdot)$ are constant, with Tamed-Euler method in (\cite{hutzenhaler_jentzen:2012}) and see in Figure \ref{SDvsTD} that for ``good'' data the two methods are close. Choosing different data, we see that Tamed-Euler (\ref{eq1000}) takes negative values, even in the first step. In particular we see, that by altering the parameters we get the results presented in Table \ref{tab:negTD} and shown in Figure \ref{negativeTD}. Note that if the Tamed-Euler takes a negative value, it explodes in the next step, because of the $3/2-$term while taking the value zero in a step results in zero terms for all the following steps.

\begin{table}[htbp] 
\centering
		\begin{tabular}{|c|c|c|}
		\hline  Set of Parameters  & Time of first & Value of \\
			$(x_0,k_1,k_2,k_3,\D,T)$ &   negative step & step \\
    \hline          $(1,1,1000,1,10^{-3},1)$ & 1 & $-0.18$\\
    \hline          $(1,1000,1,1,10^{-3},1)$ & 27 & $-17.69$ \\
		\hline
  \end{tabular}
	\caption{\small Negative values of Tamed-Euler scheme  (\ref{eq1000}) for Heston $3/2-$model.}
	\label{tab:negTD}
\end{table}

\begin{figure}[ht]
  \caption{ Difference between the Semi-Discrete scheme and Tamed-Euler scheme (\ref{eq1000}) for $x_0 = 1, k_1=1, k_2=4, k_3=1, \D=10^{-3}, T=1$.}
  \centering
    \includegraphics[width=0.5\textwidth]{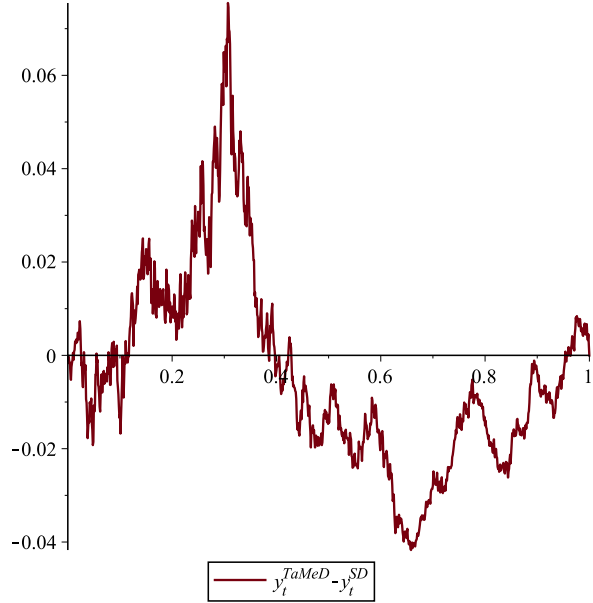}
\label{SDvsTD}
\end{figure}

\begin{figure}[ht]
 \caption{ Tamed-Euler method (\ref{eq1000}) does not preserve positivity,  $x_0=1, k_1=1000, k_2=4, k_3 = 1, \D=10^{-3}, T=1$.}
  \centering
   \includegraphics[width=0.5\textwidth]{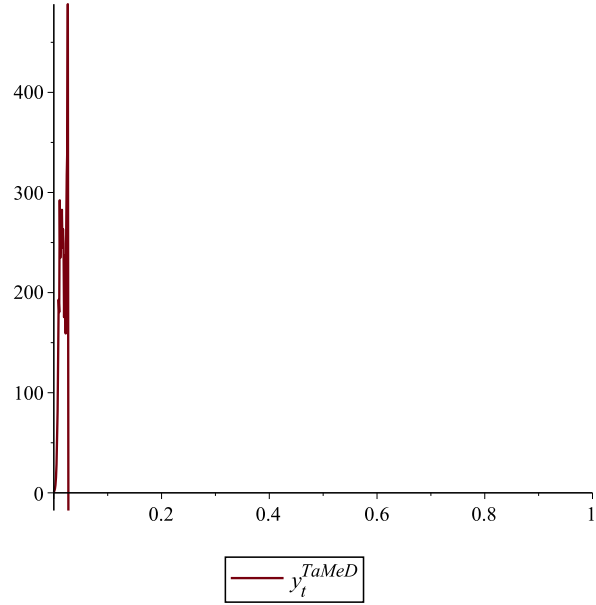}\label{negativeTD}
\end{figure}

\end{itemize}
\ere

\subsection{Example II}
Consider the following stochastic differential equation (SDE),
\beqq\label{eq31}
x_t=x_0 + \int_{0}^{t}(k_1(s) x_s -k_2(s) x_s^{2r-1})ds + \int_{0}^{t} k_3(s)x_s^{r}dW_s,\quad t\in [0,T],
\eeqq
where $x_0$ is independent of all $\{W_t\}_{0\leq t\leq T}, x_0\in \bbl^p(\Omega,\bbR)$ for some $2<p\leq\frac{r-1}{4(3-2r)} + \frac{r-1}{2(3-2r)}\frac{k_{2,min}}{(k_{3,max})^2}$ and $x_0>0,$ a.s., $k_1(\cdot), k_2(\cdot), k_3(\cdot)$ are positive and bounded functions with $2k_{2,\min}>\frac{25-9r}{r-1}k_{3,\max}^2$ and $1<r<3/2.$  

\ble\label{l31}[Positivity of $(x_t)$] In the previous setting it holds that $x_t > 0$ a.s. \ele

\bpf[Proof of Lemma \ref{l31}] Set the stopping time $\theta_R=\inf\{t\in [0,T]: x_t^{-1}>R\},$ for some $R>0,$ with the convention that $\inf\emptyset=\infty.$ Application of Ito's formula on
$x_{t\wedge\theta_R}^{-2}$ implies,
\beao
&&(x_{t\wedge \theta_R})^{-2} = (x_0)^{-2} + \int_{0}^{t\wedge \theta_R}(-2)x_s^{-3}(k_1(s)x_s- k_2(s)x_s^{2r-1}) ds\\
&&+  \int_{0}^{t\wedge\theta_R}\frac{(-2)(-3)}{2}(x_s)^{-4}k_3^2(s)x_s^{2r}ds + \int_{0}^{t\wedge \theta_R}(-2)k_3(s) (x_s)^{-3} x_{s}^{r}dW_s\\
&=& (x_0)^{-2} + \int_{0}^{t\wedge \theta_R}(-2)k_1(s)x_s^{-2}  + 2k_2(s)x_s^{2r-4} + 3k_3^2(s)x_s^{2r-4})ds\\
&&+\int_{0}^{t}(-2)k_3(s)x_s^{r-3}\,\bbi_{(0,t\wedge \theta_R)}(s)dW_s\\
&=& (x_0)^{-2} + \int_{0}^{t\wedge \theta_R}\left(-2k_1(s) x_s^{-2}  + (2k_2(s) + 3k_3^2(s))\left( x_s^{2r-4} \bbi_{(0,1]}(x_s) + x_s^{2r-4} \bbi_{(1,\infty]}(x_s)\right)\right)ds\\
&&+ M_t\\
&\leq& (x_0)^{-2} + 2k_{2,\max}T+ 3k_{3,\max}^2T + \int_{0}^{t}(2k_2(s)
+3k_3^2(s))x_s^{-2}\,\bbi_{(0,t\wedge \theta_R)}(s)ds + M_t,
\eeao
where
$$
M_t:=\int_{0}^{t}(-2)k_3(s) x_s^{r-3}\,\bbi_{(0,t\wedge
\theta_R)}(s)dW_s.
$$
Taking expectations in the above inequality and using the fact that $\bfE M_t=0,$\footnote{The function $h(u)=(-2)k_3(u) x_u^{r-3}\,\bbi_{(0,t\wedge
\theta_R)}(u)$ belongs to the space $\bbm^2([0,t];\bbR)$ thus (\cite[Theorem 1.5.8]{mao:1997}) implies $\bfE M_t=0$.} we get that
\beao
\bfE(x_{t\wedge \theta_R}^{-2}) &\leq& \bfE(x_0)^{-2} + 2k_{2,\max}T+ 3k_{3,\max}^2T + (2k_{2,\max} +3k_{3,\max}^2)\int_{0}^{t}\bfE (x_{s\wedge \theta_R})^{-2}ds\\
&\leq& \left(\bfE(x_0)^{-2} + 2k_{2,\max}T+ 3k_{3,\max}^2 T
\right)e^{(2k_2 +3k_3^2)T}<C, \eeao
where we have used Gronwall inequality with $C$ independent of $R.$ We
have that
\beqq\label{eq38}
(x_{t\wedge\theta_R})^{-2}=(x_{\theta_R})^{-2}\bbi_{(\theta_R\leq
t)}+ (x_{t})^{-2}\bbi_{(t < \theta_R)}=R^{2}\bbi_{(\theta_R\leq
t)}+ (x_{t})^{-2}\bbi_{(t<\theta_R)}. \eeqq

By relation (\ref{eq38})  we have that,
$$
\bfE\left(\frac{1}{x_{t\wedge\theta_R}^{2}}\right)=R^2\bfP(\theta_R\leq t) +
\bfE\left(\frac{1}{x_{t}^{2}}\bbi_{(t<\theta_R)}\right)<C,
$$
thus
$$
\bfP(x_t\leq0)=\bfP\left(\bigcap_{R=1}^\infty \Big\{x_t<
\frac{1}{R} \Big\}\right)=\lim_{\Rto}\bfP\left( \Big\{x_t<
\frac{1}{R} \Big\}\right)\leq  \lim_{\Rto}\bfP(\theta_R\leq t)=0.
$$
We conclude that $x_t > 0$ a.s. \epf

The following Lemma shows uniform bounds of $p-$moments of $(x_t).$
\ble\label{l201}
In the previous setting it holds that
$$
\bfE(\sup_{0\leq t\leq T}(x_t)^p)<A_1,
$$
for some $A_1>0$ and any $2< p \leq \frac{3}{2}-r + \frac{k_{2,min}}{(k_{3,max})^2}.$ \ele 
\bpf[Proof of Lemma \ref{l201}] We follow the same lines as in the proof of Lemma \ref{l2}.  In particular, we first get the bound 
$$
J(s,x):=\frac{xa(s,x) + (p-1)b^2(s,x)/2}{1+x^2}\leq\frac{k_{1,\max}x^2 + \Big(0.5(p-1)(k_{3,\max})^2 - k_{2,\min}\Big)x^{2r}}{1+x^2}\leq k_{1,\max},
$$
where the last inequality is valid for all $p$ such that $p \leq 1 + 2k_{2,\min}/(k_{3,\max})^2$ which implies 
$$
\bfE (x_t)^p\leq 2^{(p-2)/2}(1+\bfE (x_0)^p)e^{Cpt},
$$
for any $2<p \leq 1 + 2k_{2,\min}/(k_{3,\max})^2$ and all $t\in[0,T]$.
Using Ito's formula on $(x_t)^p,$ with $p \leq \frac{3}{2}-r + \frac{k_{2,\min}}{(k_{3,\max})^2}$ (in order to use Doob's martingale inequality later) we have that
\beao
&&(x_{t})^p \leq (x_0)^p + p\int_{0}^{t}\left[k_1(s)(x_s)^{p} + \left(\frac{p-1}{2}k_{3,\max}^2K^2_\phi-k_2\right)(x_s)^{p+2r-2}\right]ds + M_t\\
&\leq& (x_0)^p + p\int_{0}^{t}k_1(s)(x_s)^{p}ds + M_t,
\eeao
where $M_t=\int_{0}^{t}pk_3(s)(x_s)^{p+2r-1}dW_s.$ Taking the supremum and then expectations in the above inequality we get
$$
\bfE(\sup_{0\leq t\leq T}(x_{t})^p) \leq\left(\bfE(x_0)^p + \sqrt{4\bfE M_T^2}\right)e^{pk_{1,\max}T}:=A_1,
$$
where in the last step we have used Doob's martingale inequality to the diffusion term $M_t$\footnote{The function $h(u)=pk_3(u)\phi (x_u)(x_u)^{p+2r-1}$ belongs to the family $\bbm^2([0,T];\bbR)$ thus (\cite[Theorem 1.5.8]{mao:1997}) implies $\bfE M_t^2=\bfE(\int_0^t h(u)dW_u)^2=\bfE \int_0^t h^2(u)du,$ i.e. $M_t\in\bbl^2(\W;\bbR)$.}  and Gronwall inequality.
\epf

Model (\ref{eq31}) has super linear drift and diffusion coefficients. We study the numerical approximation of (\ref{eq31}). We propose the following Semi-Discrete numerical scheme for the transformed process $z_t=x_t^{2r-2},$ of (\ref{eq31}),
\beqq\label{eq32}
y_t=y_n + \int_{t_n}^{t} (K_1(s) - K_2(s) y_{t_n})y_sds + \int_{t_n}^{t} K_3(s)\sqrt{y_{t_n}}y_s dW_s,\quad
t\in[t_n, t_{n+1}],
\eeqq
where $y_n=y_n(t_n),$  for $n\leq T/\D$ and $y_0=x_0,$ a.s., where
\beqq\label{eq33}
K_1(s)=(2r-2)k_1(s), \,\,K_2(s)=(2r-2)k_2(s) - \frac{(2r-2)(2r-3)}{2}k_3^2(s), \,\,K_3(s)=(2r-2)k_3(s),
\eeqq
or in a  more compact form,
\beqq\label{eq34}
y_t=y_0 + \int_{0}^{t}(K_1(s) - K_2(s) y_{\hat{s}})y_sds +
\int_{0}^{t}K_3(s)\sqrt{y_{\hat{s}}}y_s dW_s,
\eeqq
where $\hat{s}=t_n,$ when $s\in[t_n,t_{n+1}).$ The linear SDE (\ref{eq34}) has a solution which, by use of Ito's formula, has the explicit form
\beqq\label{eq35}
y_t=x_0\exp\Big\{\int_{0}^{t}\left(K_1(s) - K_2(s) y_{\hat{s}} - K_3^2(s)\frac{y_{\hat{s}}}{2}\right) ds +
\int_{0}^{t}K_3(s)\sqrt{y_{\hat{s}}}dW_s\Big\},
\eeqq
where $y_t=y_t(t_0,x_0).$

\textbf{The transformation of (\ref{eq31}).} Application of Ito's formula to the function $z(t,x)=x^{2r-2},$ implies
\beam
\nonumber
z_t &=& z_0 + \int_{0}^{t}\left[(2r-2)x_s^{2r-3}(k_1(s)x_s- k_2(s) x_s^{2r-1}) + \frac{(2r-2)(2r-3)}{2}x_s^{2r-4}k_3^2(s)x_s^{2r}\right]ds\\ \nonumber
&&+ \int_{0}^{t}(2r-2)k_3(s) x_s^{2r-3}x_{s}^{r}dW_s\\ \nonumber
&=& z_0 + \int_{0}^{t}\left[k_1(s)(2r-2)x_s^{2r-2}- (2r-2)k_2(s) x_s^{4r-4} + \frac{(2r-2)(2r-3)}{2}k_3^2(s)x_s^{4r-4}\right]ds\\
\label{eq37}&&+ \int_{0}^{t}(2r-2)k_3(s) x_s^{3r-3}dW_s\\
\nonumber&=& z_0 + \int_{0}^{t}(K_1(s)z_s - K_2(s)z_s^2)ds +\int_{0}^{t}K_3(s) z_s^{3/2}dW_s,
\eeam
where $K_1(\cdot), K_2(\cdot), K_3(\cdot)$ are given by (\ref{eq33}).

In order to use Proposition \ref{pr0} we have to verify that
$$
K_1(s)>0, \quad K_2(s)>0, \quad K_3(s)>0, \quad 2K_{2,\min}>7K_{3,\max}^2.
$$
Since $1< r < 3/2$ we immediately have $K_1(s)>0$ and $K_3(s)>0.$ Moreover
$$
K_2(s)=(2r-2)k_2(s) - \frac{(2r-2)(2r-3)}{2}k_3^2(s)>\frac{(2r-2)}{2}k_{3,\max}^2(4 -2r)>0,
$$
and is easy to see that \beao 2K_{2,\min}> 7K_{3,\max}^2. \eeao

\bpr \label{pr31} In the previous setting, the following convergence to the true solution
of (\ref{eq31}) in the mean square sense holds, \beqq \label{eq36}
\lim_{\D\rightarrow0}\bfE\sup_{0\leq t\leq
T}|y_t^{\frac{1}{2r-2}}-x_t|^2=0.\eeqq\epr

\subsubsection*{Proof of Proposition \ref{pr31}}
In order to prove Proposition \ref{pr31} we first transform the
original SDE (\ref{eq31}) to a SDE (\ref{eq11}), later on verify
the assumptions of Example I to use Proposition \ref{pr0}, and in
the end make the necessary arrangements for the approximation of
the original SDE.

\subsubsection{Convergence result}
We use the following inequality implied by the mean value theorem
$$
|y_t^{\frac{1}{2r-2}}-x_t| = |y_t^{\frac{1}{2r-2}} - z_t^{\frac{1}{2r-2}}|\leq \frac{1}{2r-2}\left( |y_t|^{\frac{1}{2r-2}-1} + |z_t|^{\frac{1}{2r-2}-1} \right) |z_t-y_t|.
$$
thus we get that
$$
|y_t^{\frac{1}{2r-2}}-x_t|^2\leq\frac{2}{(2r-2)^2} \left( |y_t|^{\frac{3-2r}{r-1}} + |z_t|^{\frac{3-2r}{r-1}} \right) |z_t-y_t|^2.
$$
Set the stopping time $\theta_R=\inf\{t\in [0,T]: |y_t|>R \, \mbox{or} \, |x_t|> R\},$ for some $R>0$ big enough. Taking the supremum and then expectations in the above inequality yields,
\beao
&&\bfE\sup_{0\leq t\leq T}|y_t^{\frac{1}{2r-2}}-x_t|^2 \leq c_r\Big[\bfE\sup_{0\leq t\leq T}\left( |y_{t\wedge\theta_R}|^{\frac{3-2r}{r-1}} + |z_{t\wedge\theta_R}|^{\frac{3-2r}{r-1}} \right) |z_{t\wedge\theta_R}-y_{t\wedge\theta_R}|^2\\
&&+
\bfE\sup_{0\leq t\leq T} \left( |y_t|^{\frac{3-2r}{r-1}} + |z_t|^{\frac{3-2r}{r-1}} \right) |z_t-y_t|^2\bbi_{(\theta_R\leq t)}\Big]\\
&\leq& c_{r,R}\bfE\sup_{0\leq t\leq T }|z_{t\wedge\theta_R}-y_{t\wedge\theta_R}|^2 +
c_r\frac{2\del}{p}\bfE\sup_{0\leq t\leq T}\left( |y_t|^{\frac{3-2r}{r-1}} + |z_t|^{\frac{3-2r}{r-1}} \right)^{p/2} |z_t - y_t|^{p}\\
&& + c_r\frac{(p-2)}{p\del^{2/(p-2)}}\bfP(\theta_R\leq T),
\eeao
where in the second step we have applied Young inequality,
$$
ab \leq \frac{\del}{w}a^w + \frac{1}{q\del^{q/w}}b^q,
$$
for $a=\sup_{0\leq t\leq T}\left( |y_t|^{\frac{3-2r}{r-1}} +
|z_t|^{\frac{3-2r}{r-1}} \right) |z_t-y_t|^2,
b=\bbi_{(\theta_R\leq t)}, w=p/2, q=p/(p-2), \del>0,$ and
$$
c_r=\frac{2}{(2r-2)^2}, \quad c_{r,R}=2c_rR^{\frac{3-2r}{r-1}}.\footnote{For all $t<\theta_R$ it holds that $|x_t|\leq R$ or $|z_t|\leq R$.}
$$
It holds that
$$
\!\!\bfP(\theta_R\leq T)\!\leq\!\bfE\!\left(\!\bbi_{(\theta_R\leq T)}\frac{|y_{\theta_R}|^{p}}{R^p}\right)\!+\!\bfE\left(\!\bbi_{(\theta_R\leq T)}\frac{|x_{\theta_R}|^{p}}{R^p}\right)\!\leq\!\frac{1}{R^p}\!\left(\!\bfE\sup_{0\leq t\leq T}|y_t|^{p}\!+\!\bfE\sup_{0\leq t\leq T}|x_t|^{p}\right)\!\leq\!\!\frac{2A}{R^p},
$$
where $A$ is the maximum of the bounding moment constants of $y$ and $x.$
Moreover, we have that,
\beao
&&\bfE\sup_{0\leq t\leq T}\left( |y_t|^{\frac{3-2r}{r-1}} + |z_t|^{\frac{3-2r}{r-1}} \right)^{p/2} |z_t - y_t|^{p}\leq 2^{\frac{3p}{2}-2} \bfE\sup_{0\leq t\leq T}\left( |y_t|^{\frac{(3-2r)p}{2(r-1)}} + |z_t|^{\frac{(3-2r)p}{2(r-1)}} \right) (|z_t|^{p} + |y_t|^{p})\\
&\leq& 2^{\frac{3p}{2}-2} \bfE\sup_{0\leq t\leq T}\left(|y_t|^{\frac{(3-2r)p}{2(r-1)}}|z_t|^{p} + |y_t|^{(\frac{3-2r}{2(r-1)}+1)p}  + |z_t|^{\frac{(3-2r)p}{2(r-1)}}|y_t|^{p} + |z_t|^{(\frac{3-2r}{2(r-1)}+1)p}\right)\\
&\leq& 2^{\frac{3p}{2}-2} \bfE\sup_{0\leq t\leq T}\left(\frac{|y_t|^{\frac{3-2r}{r-1}p}}{2} +\frac{|z_t|^{2p}}{2} + |y_t|^{\frac{p}{2(r-1)}} +  \frac{|z_t|^{\frac{3-2r}{r-1}p}}{2} + \frac{|y_t|^{2p}}{2}  + |z_t|^{\frac{p}{2(r-1)}}\right),
\eeao
where we have used again Young inequality. When $\frac{5}{4}<r<\frac{3}{2}$ we have that 
$\frac{3-2r}{r-1}<\frac{1}{2(r-1)}<2,$ thus it suffices to bound the moments of $|z_t|^{2p}$ and $|y_t|^{2p}.$ Note that by Lemma \ref{l2} the uniform bound for the moment of $(z_t)^{2p}$ holds when $2< p \leq \frac{k_{2,min}}{2(k_{3,max})^2}$ and by  Lemma \ref{l3} the uniform bound for the moment of $(y_t)^{2p}$ is valid for any $2<p \leq \frac{1}{8} + \frac{k_{2,min}}{4(k_{3,max})^2},$ thus for $2< p \leq \frac{k_{2,min}}{2(k_{3,max})^2}\bigwedge \frac{1}{8} + \frac{k_{2,min}}{4(k_{3,max})^2}$\footnote{\label{f1}We also have to ensure that Lemma \ref{l201} holds, thus we have to choose $p$ such that $2<p\leq\frac{3}{2}-r + \frac{k_{2,min}}{(k_{3,max})^2}\bigwedge\frac{k_{2,min}}{2(k_{3,max})^2}\bigwedge\frac{1}{8} + \frac{k_{2,min}}{4(k_{3,max})^2}$ or equivalently we have to choose $p$ such that $2<p\leq\frac{1}{8} + \frac{k_{2,min}}{4(k_{3,max})^2}$ whose existence is ensured by the condition $2k_{2,min}\geq15(k_{3,max})^2.$} we get that $\bfE\sup_{0\leq t\leq T}\left(|z_t|^{2p} +  |y_t|^{2p}\right)<A,$ for some $A>0.$ In the case $1<r<\frac{5}{4}$ it suffices to bound the moments of $|z_t|^{\frac{3-2r}{r-1}p}$ and $|y_t|^{\frac{3-2r}{r-1}p}.$ Again by Lemma \ref{l2} the uniform bound for the moment of $|z_t|^{\frac{3-2r}{r-1}p}$ holds when $2< p \leq \frac{r-1}{3-2r}\frac{k_{2,min}}{(k_{3,max})^2}$ and by  Lemma \ref{l3} the uniform bound for the moment of $|y_t|^{\frac{3-2r}{r-1}p}$ is valid for any $2<p \leq \frac{r-1}{4(3-2r)} + \frac{r-1}{2(3-2r)}\frac{k_{2,min}}{(k_{3,max})^2},$ thus for $2< p \leq \frac{r-1}{3-2r}\frac{k_{2,min}}{(k_{3,max})^2}\bigwedge \frac{r-1}{4(3-2r)} + \frac{r-1}{2(3-2)r}\frac{k_{2,min}}{(k_{3,max})^2}$\footnote{\label{f2}We also have to ensure that Lemma \ref{l201} holds, thus we have to choose $p$ such that $2<p\leq\frac{3}{2}-r + \frac{k_{2,min}}{(k_{3,max})^2}\bigwedge\frac{r-1}{3-2r}\frac{k_{2,min}}{(k_{3,max})^2}\bigwedge\frac{r-1}{4(3-2r)} + \frac{r-1}{2(3-2r)}\frac{k_{2,min}}{(k_{3,max})^2}$ or equivalently we have to choose $p$ such that $2<p\leq\frac{r-1}{4(3-2r)} + \frac{r-1}{2(3-2r)}\frac{k_{2,min}}{(k_{3,max})^2}$ whose existence is ensured by the condition $2k_{2,min}\geq\frac{25-9r}{r-1}(k_{3,max})^2.$}
we get that $\bfE\sup_{0\leq t\leq T}\left(|z_t|^{\frac{3-2r}{r-1}p} +  |y_t|^{\frac{3-2r}{r-1}p}\right)<A,$ for some $A>0.$ Thus, by Footnotes \ref{f1} and \ref{f2} and the condition 
$2k_{2,min}\geq\left(\frac{25-9r}{r-1} \bigvee 15\right)(k_{3,max})^2$ or equivalently 
$2k_{2,min}\geq\frac{25-9r}{r-1}(k_{3,max})^2$ we get the bound $\bfE\sup_{0\leq t\leq T}\left( |y_t|^{\frac{3-2r}{r-1}} + |z_t|^{\frac{3-2r}{r-1}} \right)^{p/2} |z_t - y_t|^{p}<C(p)A,$ where
$C(p)$ is a constant depending on $p.$ Collecting all the estimates together,
\beao
\bfE\sup_{0\leq t\leq T}|y_t^{\frac{1}{2r-2}}-x_t|^2 &\leq& c_{r,R}\bfE\sup_{0\leq t\leq T }|z_{t\wedge\theta_R}-y_{t\wedge\theta_R}|^2 + c_r\frac{C(p)A}{p}\del + c_r\frac{2(p-2)A}{p}\frac{1}{\del^{2/(p-2)}R^p}\\
&&:=I_1 + I_2 + I_3.
\eeao
Given any $\ep>0,$ we may first choose $\del$ such that $I_2<\ep/3,$ then choose $R$ such that  $I_3<\ep/3,$ and finally $\D$ such that $I_1<\ep/3,$ which is justified  by Proposition \ref{pr0} to get that  $\bfE\sup_{0\leq t\leq T}|y_t^{\frac{1}{2r-2}}-x_t|^2 <\ep,$ as required to verify (\ref{eq36}).

\bre
Proposition \ref{pr31} implies that
our explicit numerical scheme converges in the mean square sense.
Moreover, we get that our numerical scheme preserves
positivity. Example (\ref{eq31}) covers  super-linear problems both in drift and diffusion.
\ere

\subsection{Example III}
Consider the following stochastic differential equation (SDE),
\beqq\label{eq41} x_t=x_0 + \int_{0}^{t}(k_1(s) x_s -k_2(s)
x_s^{q})ds + \int_{0}^{t} k_3(s)x_s^{r}\phi(x_s)dW_s,\quad t\in
[0,T], \eeqq where $\phi(\cdot)$ is a locally Lipschitz and
bounded function with locally Lipschitz constant $C_R^\phi,$
bounding constant $K_\phi, x_0$ is independent of all
$\{W_t\}_{0\leq t\leq T}, x_0\in \bbl^p(\Omega,\bbR)$ for every
$2<p, \bfE |\ln x_0|<\infty$ and $x_0>0,$ a.s., $k_1(\cdot),
k_2(\cdot), k_3(\cdot)$ are positive and bounded functions and $q$
is odd with $q>2r-1$ where $3/2<r<2.$  The above conditions on the
parameters imply the uniform bound of $|x_t|^p$ as shown in the following result.

\ble\label{l42}[Moment bound for original SDE]
In the previous setting it holds that
$$
\bfE(\sup_{0\leq t\leq T}|x_t|^p)<A_1,
$$
for some $A_1>0$ and every $p > 2.$  \ele
\bpf[Proof of Lemma \ref{l42}] In the case of $x$'s outside a finite ball of radius $R,$ with $R>1,$ and when $s\in[0,T]$ we have that
\beao
&&J(s,x):=\frac{xa(s,x) + (p-1)b^2(s,x)/2}{1+x^2}=\frac{x(k_1(s) x-k_2(s) x^{q}) + (p-1)k_3^2(s)[x^{r}\phi(x)]^2/2}{1+x^2}\\
&=&\frac{k_1(s)x^2 - k_2(s)x^{q+1} + 0.5(p-1)k_3^2(s)x^{2r}\phi^2(x)}{1+x^2}\\
&\leq& k_{1,\max}, \eeao
where the the last inequality is valid for all $p>2$ and we have used $q+1>2r$ and that $q$ is odd. Thus $J(s,x)$ is bounded for all $(s,x)\in[0,T]\times\bbR,$ since when $|x|\leq R$ we have that $J(s,x)$ is finite and say $J(s,x)\leq
C.$ Application of (\cite[Theorem 2.4.1]{mao:1997}) implies
$$
\bfE |x_t|^p\leq 2^{(p-2)/2}(1+\bfE |x_0|^p)e^{C pt},
$$
for any $2<p$ and all $t\in[0,T]$. Using Ito's formula on $|x_t|^p,$
we have that
\beao
&&|x_{t}|^p = |x_0|^p + \int_{0}^{t}p|x_s|^{p-2}x_s(k_1(s)x_s- k_2(s) x_s^{q}) ds\\
&&+  \int_{0}^{t}\frac{p}{2}\left(|x_s|^{p-2} + (p-2)|x_s|^{p-4}x_s^2\right)[k_3(s)x_s^{r}\phi (x_s)]^2 ds + \int_{0}^{t}pk_3(s)|x_s|^{p-2}x_sx_s^{r}\phi (x_s)dW_s\\
&\leq& |x_0|^p + p\int_{0}^{t}\left[k_1(s) - k_2(s)(x_s)^{q-1} +\frac{p-1}{2}k_3^2(s)K^2_\phi (x_s)^{2r-2} \right]|x_s|^{p}ds \\
&&+ \int_{0}^{t}pk_3(s)\phi (x_s)|x_s|^p(x_s)^{r-1}dW_s\\
&\leq& |x_0|^p + C\int_{0}^{t} |x_s|^{p}ds + M_t,
\eeao
where we have used that $0<2r-2<q-1, q$ is odd and $M_t=\int_{0}^{t}pk_3(s)\phi (x_s)|x_s|^{p}(x_s)^{r-1}dW_s.$ Taking the supremum and then expectations in the above inequality we get
\beao
\bfE(\sup_{0\leq t\leq T}|x_{t}|^p) &\leq& \bfE|x_0|^p + C\bfE\left(\sup_{0\leq t\leq T}\int_{0}^{t} |x_s|^{p}ds\right) + \bfE\sup_{0\leq t\leq T}M_t\\
&\leq&\bfE|x_0|^p + C\int_{0}^{t} \bfE(\sup_{0\leq l\leq s}|x_l|^{p})ds + \sqrt{\bfE\sup_{0\leq t\leq T}M_t^2}\\
&\leq&\left(\bfE|x_0|^p + \sqrt{4\bfE M_T^2}\right)e^{CT}:=A_1,
\eeao
where in the last step we have used Doob's martingale inequality to the diffusion term $M_t$\footnote{The function $h(u)=pk_3(u)\phi (x_u)|x_u|^{p}x_s^{r-1}$ belongs to the family $\bbm^2([0,T];\bbR)$ thus (\cite[Theorem 1.5.8]{mao:1997}) implies $\bfE M_t^2=\bfE(\int_0^t h(u)dW_u)^2=\bfE \int_0^t h^2(u)du,$ i.e. $M_t\in\bbl^2(\W;\bbR)$.}  and Gronwall inequality.
\epf

Model (\ref{eq41}) has super linear drift and diffusion coefficients. We study the numerical approximation of (\ref{eq41}).
We propose the following Semi-Discrete numerical scheme for (\ref{eq41})
\beqq\label{eq42}
y_t=y_n + \int_{t_n}^{t} (k_1(s) - k_2(s) y_{t_n}^{q-1})y_sds   + \int_{t_n}^{t} k_3(s)y_{t_n}^{r-1}\phi(y_{t_n})y_s dW_s,\quad t\in[t_n, t_{n+1}],
\eeqq
where $y_n=y_n(t_n),$  for $n\leq T/\D$ and $y_0=x_0,$ a.s.,
or in a  more compact form,
\beqq \label{eq44}
y_t=y_0 + \int_{0}^{t}(k_1(s) - k_2(s) y_{\hat{s}}^{q-1})y_sds +
\int_{0}^{t}k_3(s)y_{\hat{s}}^{r-1}\phi(y_{\hat{s}})y_s dW_s,
\eeqq
where $\hat{s}=t_n,$ when $s\in[t_n,t_{n+1}).$ The linear SDE (\ref{eq44}) has a solution which, by use of Ito's formula, has the explicit form (\cite[Chapter 4.4, relation(4.10)]{kloeden_platen:1995})
\beqq\label{eq45}
y_t=x_0\exp\Big\{\int_{0}^{t}\left(k_1(s) - k_2(s) y_{\hat{s}}^{q-1} - k_3^2(s)\frac{y_{\hat{s}}^{2r-2}\phi^2(y_{\hat{s}})}{2}\right) ds +
\int_{0}^{t}k_3(s)y_{\hat{s}}^{r-1}\phi(y_{\hat{s}})dW_s\Big\},
\eeqq
where $y_t=y_t(t_0,x_0).$
\bpr \label{pr41}
The following convergence to the true solution of
(\ref{eq41}) in the mean square sense holds,
\beqq \label{eq46}
\lim_{\D\rightarrow0}\bfE\sup_{0\leq t\leq T}|y_t-x_t|^2=0.
\eeqq
\epr

\subsubsection{Proof of Proposition \ref{pr41}}\label{ssec:s42}

In order to prove Proposition \ref{pr41} we just need to verify the assumptions of Theorem \ref{t1}. Let
\beao
a(s,x)=k_1(s) x -k_2(s) x^{q}, && f(s,r,x,y)=(k_1(s)-k_2(s) x^{q-1})y\\
b(s,x)= k_3(s)x^{r}\phi(x), && g(s,r,x,y)=k_3(s)x^{r-1}\phi(x)y.
\eeao
We verify Assumption A for $f.$ The conditions on the parameters imply that $q>2.$
Let $R>0$ such that $|x_1|\vee|x_2|\vee|y_1|\vee|y_2|\vee|s|\vee|r|\leq R.$ We have that
\beao
&&|f(s,r,x_1,y_1)-f(s,r,x_2,y_2)|=|(k_1(s)-k_2(s) x_1^{q-1})y_1 - (k_1(s)-k_2(s) x_2^{q-1})y_2|\\
&\leq& |k_1(s)||y_1-y_2| + |k_{2,\max}|(|x_2|^{q-1}|y_1-y_2| + |y_1||x_1^{q-1}-x_2^{q-1}|)\\
&\leq& (|k_{1,\max}| + |k_{2,\max}|R^{q-1})|y_1-y_2| + |k_{2,\max}|R|x_1^{q-1}-x_2^{q-1}|\\
&\leq& (|k_{1,\max}| + |k_{2,\max}|R^{q-1})|y_1-y_2| + 2|k_{2,\max}| (q-1)R^{q-1}|x_1-x_2|\\
&\leq& C_R\left( |x_1-x_2| + |y_1-y_2| \right),
\eeao
where we have applied the mean value theorem for the function $x^{q-1}$, thus Assumption A holds for $f$ with $C_R:=(|k_{1,\max}| + |k_{2,\max}|R^{q-1})\vee(2|k_{2,\max}|(q-1)R^{q-1}).$

We verify Assumption A for $g.$ Since $1/2<r-1<1$ we have that $g_1(x)=x^{r-1}$ is locally
$1/2-$Holder continuous in $x,$ i.e. \beqq\label{eq47} |g_1(x_1)-
g_1(x_2)| \leq C_R \sqrt{|x_1 - x_2|}. \eeqq Let $R>0$ such that
$|x_1|\vee|x_2|\vee|y_1|\vee|y_2|\vee|s|\vee|r|\leq R.$ We have
that \beao
&&|g(s,r,x_1,y_1)-g(s,r,x_2,y_2)|=\left|k_3(s) x_1^{r-1}\phi (x_1)y_1 -k_3(s) x_2^{r-1}\phi(x_2)y_2\right|\\
&\leq&|k_{3,\max}| \left(|x_1|^{r-1}| \phi (x_1)| |y_1-y_2| + |y_2| \big| x_1^{r-1}\phi(x_1) - x_1^{r-1}\phi( x_2) + x_1^{r-1}\phi(x_2) - x_2^{r-1}\phi(x_2)\big|\right)\\
&\leq& |k_{3,\max}|\left(K_\phi R^{r-1} |y_1-y_2| + R|x_1|^{r-1} | \phi(x_1) - \phi(x_2)| +  R K_\phi|x_1^{r-1} - x_2^{r-1}|\right)\\
&\leq& |k_{3,\max}|\left( K_\phi R^{r-1} |y_1-y_2| + R^{r} C^\phi_R| x_1 - x_2| +  R K_\phi \sqrt{|x_1 - x_2|}\right)\\
&\leq& C_R \left(|x_1 -x_2| +|y_1-y_2| + \sqrt{|x_1 -x_2|}\right),
\eeao
where we have used (\ref{eq47}) and $C_R:=|k_{3,\max}|\left(C_R^\phi R^{r}\vee
K_\phi R^{r-1}\vee K_\phi R\right).$ Thus, Assumption A holds for $g$.

\ble\label{l41}[Positivity of $(x_t)$] In the previous setting it holds that $x_t > 0$
a.s. \ele

\bpf[Proof of Lemma \ref{l41}] Set the stopping time $\theta_R=\inf\{t\in [0,T]: x_t^{-1}>R\},$ for some $R>0,$ with the convention that $\inf\emptyset=\infty.$ Application of Ito's formula on
$\ln x_{t\wedge\theta_R}$ implies,
\beao
&&\ln x_{t\wedge \theta_R} = \ln x_0 + \int_{0}^{t\wedge \theta_R}\frac{1}{x_s}(k_1(s)x_s- k_2(s) x_s^{q}) ds +\int_{0}^{t\wedge\theta_R}\left(-\frac{1}{x_s^2}\right)k_3^2(s)x_s^{2r}\phi^2(x_s)ds\\
&&+ \int_{0}^{t\wedge \theta_R}\frac{1}{x_s}k_3(s) x_{s}^{r}\phi(x_s)dW_s\\
&=& \ln x_0 + \int_{0}^{t\wedge \theta_R}\left(k_1(s) - k_2(s)x_s^{q-1} - k_3^2(s)x_s^{2r-2}\phi^2(x_s)\right)ds+\int_{0}^{t\wedge \theta_R}k_3(s)x_s^{r-1}\phi(x_s)dW_s.
\eeao
Taking absolute values in the above equality and then expectations and using Jensen inequality and then Ito's isometry on the diffusion term, $M_t,$ we get
\beao
&&\bfE |\ln x_{t\wedge \theta_R}|\leq \bfE |\ln x_0| + T(|k_{1,\max}| + |k_{2,\max}|\bfE\sup_{0\leq t\leq T}|x_t|^{q-1} +  |k_{3,\max}|^2K_\phi^2\bfE\sup_{0\leq t\leq T}|x_t|^{2r-2})\\
&& + \bfE|M_t|\\
&\leq& \bfE |\ln x_0| + (|k_{1,\max}| + (|k_{2,\max}|+|k_{3,\max}|^2)A_1 + |k_{3,\max}|^2K_\phi^2)T + \sqrt{4\bfE M_T^2}<C, \eeao
where $A_1$ is as in Lemma \ref{l42} and $M_t:=\int_{0}^{t}k_3(s)x_s^{r-1}\phi(x_s)\,\bbi_{(0,t\wedge\theta_R)}(s)dW_s.$    Now we proceed as in Lemmata \ref{l4} and \ref{l31}, to get first that $\lim_{\Rto}\bfP(\theta_R\leq t)=0$ and then conclude that $\bfP (x_t\leq0),$ i.e. $x_t > 0$ a.s.
\epf

\ble\label{l43}[Moment bound for Semi-Discrete approximation]
In the previous setting it holds that
$$\bfE(\sup_{0\leq t\leq T}(y_t)^p)<A_2,
$$
for some $A_2>0$ and for every $p>2.$   \ele
\bpf[Proof of Lemma \ref{l43}] Set the stopping time $\theta_R=\inf\{t\in [0,T]: y_t>R \},$ for some $R>0,$ with the convention that $\inf\emptyset=\infty.$ Application of Ito's
formula on $(y_{t\wedge \theta_R})^p,$ implies,
\beao
&&(y_{t\wedge \theta_R})^p = (y_0)^p + \int_{0}^{t\wedge \theta_R}p(y_s)^{p-1}(k_1(s)- k_2(s) y_{\hat{s}}^{q-1})y_s ds\\
&&+  \int_{0}^{t\wedge\theta_R}\frac{p(p-1)}{2}(y_s)^{p-2}\left[ k_3(s)y_{\hat{s}}^{r-1}\phi (y_{\hat{s}})y_s\right]^2 ds + \int_{0}^{t\wedge \theta_R}p k_3(s)(y_s)^{p-1}y_{\hat{s}}^{r-1}\phi (y_{\hat{s}})y_sdW_s\\
&=&(x_0)^p + \int_{0}^{t\wedge \theta_R}\left(p(k_1(s)- k_2(s) y_{\hat{s}}^{q-1}) + \frac{p(p-1)k_3^2(s)}{2}y_{\hat{s}}^{2r-2}\phi^2 (y_{\hat{s}})\right)(y_s)^{p}ds\\
&& + \int_{0}^{t\wedge \theta_R}p k_3(s)y_{\hat{s}}^{r-1}\phi (y_{\hat{s}}) (y_s)^{p}dW_s\\
&\leq& (x_0)^p + p\int_{0}^{t}\left[-k_2(s)(y_{\hat{s}})^{q-1} +  \frac{p-1}{2}k_{3,\max}^2K^2_\phi y_{\hat{s}}^{2r-2} + k_{1,\max}\right](y_s)^{p}\,\bbi_{(0,t\wedge \theta_R)}(s)ds + M_t\\
&\leq& (x_0)^p + C\int_0^t(y_s)^{p}\,\bbi_{(0,t\wedge \theta_R)}(s)ds + M_t, \eeao
where we have used that $q-1>2r-2>1,$ the last inequality is valid for $p>2,$ the constant $C$ is independent of $R$  and $M_t:=\int_{0}^{t\wedge \theta_R}p k_3(s)y_{\hat{s}}^{r-1}\phi(y_{\hat{s}}) (y_s)^{p}dW_s.$ Taking expectations and using that $\bfE M_t=0$ we get
\beao
\bfE(y_{t\wedge \theta_R})^{p} &\leq& \bfE(x_0)^{p} + C\int_0^t\bfE(y_{s\wedge \theta_R})^{p}ds\\
&\leq& \bfE(x_0)^{p}e^{CT},
\eeao
where in the second step we have applied Gronwall inequality.
We have that
$$
(y_{t\wedge \theta_R})^{p}=(y_{\theta_R})^{p}\bbi_{(\theta_R\leq
t)}+ (y_{t})^{p}\bbi_{(t < \theta_R)}=R^{p}\bbi_{(\theta_R\leq
t)}+ (y_{t})^{p}\bbi_{(t<\theta_R)},
$$
thus taking expectations in the above inequality and using the
estimated upper bound for $\bfE(y_{t\wedge \theta_R})^{p}$ we
arrive at
$$
\bfE (y_{t})^{p}\bbi_{(t<\theta_R)}\leq \bfE(x_0)^{p}e^{CT}
$$
and taking limits in both sides as $\Rto$ we get that
$$
\lim_{\Rto}\bfE (y_{t})^{p}\bbi_{(t<\theta_R)}\leq \bfE(x_0)^{p}e^{CT}.
$$
Fix $t$. The sequence $(y_{t})^{p}\bbi_{(t<\theta_R)}$ is
nondecreasing in $R$ since $\theta_R$ is increasing in $R$ and
$t\wedge \theta_R\rightarrow t$ as $\Rto$ and
$(y_{t})^{p}\bbi_{(t<\theta_R)}\rightarrow(y_{t})^{p}$ as $\Rto,$
thus the monotone convergence theorem
implies
\beqq \label{eq166}
\bfE (y_{t})^{p}\leq \bfE(x_0)^{p}e^{CT}, \eeqq
for any $2<p.$ Following the same lines as in Lemma \ref{l42}, i.e. using again Ito's formula on $(y_t)^p$, taking the
supremum and then using Doob's martingale inequality on the
diffusion term we obtain the desired result. \epf

\section{Numerical Experiments.}\label{sec:s5}
\setcounter{equation}{0}

We study the numerical approximation of the following SDE,
\beqq\label{eq00001} 
x_t=x_0 + \int_{0}^{t}(k_1 x_s -k_2x_s^{2})ds + \int_{0}^{t} k_3x_s^{3/2} dW_s,\quad t\in [0,T], \eeqq 
where $x_0$ is independent of all $\{W_t\}_{0\leq t\leq T}, x_0\in \bbl^{4p}(\Omega,\bbR)$ for some $2<p$ and $x_0>0,$ a.s., $\mathbb{E}(x_0)^{-2} < A$, $k_1, k_2, k_3$ are positive constants with $k_{2}> \frac{7}{2}(k_{3})^2.$  Model (\ref{eq00001}) has super linear drift and diffusion coefficients.
 
In Proposition \ref{pr0} we have shown that the following Semi-Discrete numerical scheme\footnote{The existence and uniqueness of $y_t^{SD}$ is shown in Appendix \ref{ap1}.} (in a more general setting with time-varying coefficients) 
\beqq\label{eq00002}
y_t^{SD}=y_n + \int_{t_n}^{t} (k_1 - k_2y_{t_n})y_sds   + \int_{t_n}^{t} k_3\sqrt{y_{t_n}}y_s  dW_s,\quad t\in[t_n, t_{n+1}],
\eeqq
where $y_n=y_n(t_n),$  for $n\leq T/\D$ and $y_0=x_0,$ a.s., or in a  more compact form,
\beqq\label{eq00003}
y_t^{SD}=y_0 + \int_{0}^{t}(k_1 - k_2y_{\hat{s}})y_sds +
\int_{0}^{t}k_3\sqrt{y_{\hat{s}}}y_s dW_s,
\eeqq
where $\hat{s}=t_n,$ when $s\in[t_n,t_{n+1}),$ converges to the true solution of (\ref{eq00001}) in the mean square sense, that is
\beqq \label{eq00004}
\lim_{\D\rightarrow0}\bfE\sup_{0\leq t\leq T}|y_t^{SD}-x_t|^2=0.
\eeqq
  
Relation (\ref{eq00004}) does not show the order of convergence. We aim to show 
experimentally the order.  

The linear SDE (\ref{eq00003}) has a solution which, by use of Ito's formula, has the explicit form
\beqq\label{eq00005}
y_t^{SD}=x_0\exp\Big\{\int_{0}^{t}\left(k_1 - k_2y_{\hat{s}} - k_3^2\frac{y_{\hat{s}}}{2}\right) ds +
\int_{0}^{t}k_3\sqrt{y_{\hat{s}}} dW_s\Big\},
\eeqq
where $y_t=y_t(t_0,x_0).$ The Semi-Discrete numerical scheme preserves positivity, which is a desirable modeling property. 

In order to estimate the endpoint error $\ep=\bfE{|y_T-x_T|},$ where $x_T$ is the exact solution of (\ref{eq00001}) and $y_T$ is the Semi-Discrete approximation (\ref{eq00005}) we follow a standard procedure (\cite[Section 3.3]{kloeden_platen_schurz:2003}). We compute $M$ batches of $L$ simulation paths. Each batch is estimated by 
$$
\hat{\ep_j}=\frac{1}{L}\sum_{i=1}^L|y_T^{i,j} - x_T^{i,j}|
$$
and the Monte Carlo estimator of the error
$$
\hat{\ep}=\frac{1}{M}\sum_{j=1}^M \hat{\ep_j}=\frac{1}{ML}\sum_{j=1}^M\sum_{i=1}^L|y_T^{i,j} - x_T^{i,j}|,
$$
requires $M\cdot L$ Monte Carlo sample paths. When the batch size averages $L\geq15$ they can be considered as Gaussian. A $100(1-\alpha)\%$ confidence interval for the error $\ep$ is of the form
$$
\left(\hat{\ep} - t_{1-\alpha,M-1}\cdot\sqrt{\frac{1}{M(M-1)}\sum_{j=1}^{M}(\hat{\ep}_j-\hat{\ep})^2}, \hat{\ep} + t_{1-\alpha,M-1}\cdot\sqrt{\frac{1}{M(M-1)}\sum_{j=1}^{M}(\hat{\ep}_j-\hat{\ep})^2}\right).
$$
	We simulate $20\cdot 100=2000$ paths\footnote{\label{simu}We simulate with $3.06$GHz Intel Pentium, $1.49$GB of RAM in Maple $16$ Software. The effort made is just for the purpose of the order of convergence and not for the efficiency of the computer code-time.}. The choice for $L=100$ is considered in (\cite[p.118]{kloeden_platen_schurz:2003}).  We should not forget to change the student t-test quantile $t_{1-\alpha,M-1}$ when we change the number $M$ of batches or the significance level $\alpha.$ For example for the $90\%$ confidence intervals we have 
	\begin{table}[htbp]
\centering
		\begin{tabular}{|c|c|c|c|c|c|c|c|}
		\hline t-test quantile & $M=10$ & $M=20$ & $M=30$ & $M=40$ & $M=60$ & $M=100$ & $M=200$\\	
		\hline $t_{0.9,M-1}$ & $1.83$ & $1.73$ & $1.70$ & $1.68$ & $1.67$ & $1.66$ & $M=1.65$\\
	\hline	\end{tabular}
	\caption{t-test quantiles, batches, level of confidence.}
	\label{tab:t}
\end{table}

We discretize with a number of steps in power of $2.$ The iterative SD-procedure reads
$$
y_{t_{n+1}}^{SD}=y_{t_{n}}\exp\Big\{\left(k_1 - k_2y_{t_n} - \frac{k_3^2 y_{t_n}}{2}\right)\D +
k_3\sqrt{y_{t_n}}\D W_n\Big\},
$$
for $n=0,\ldots,N-1,$ where $\D W_n:=W_{t_{n+1}}-W_{t_{n}}$ are the increments of the Brownian motion.

We want to compare our results with two other methods. The first is an implicit Milstein scheme proposed in (\cite[Section 2.2]{higham_et_al:2013}), which takes the form
\beao
y_{t_{n+1}}^{HMS}&=&\frac{1}{2(k_2 + \frac{3}{4}(k_3)^2)\D}\Big(-(1-k_1\D)\\
&&+ \sqrt{(1-k_1\D)^2 + 4(k_2 + \frac{3}{4}(k_3)^2)\D(y_{t_{n}} + k_3y_{t_{n}}^{3/2}\D W_n + \frac{3}{4}(k_3)^2y_{t_{n}}^2(\D W_n)^2} \Big)
\eeao
and the second is a tamed Euler-Maruyama scheme proposed in  (\cite[Relation 4]{hutzenhaler_jentzen:2012}), which reads 
$$
y_{t_{n+1}}^{TAMeD}=y_{t_{n}} + \frac{(k_1y_{t_{n}}- k_2y_{t_{n}}^2)\D + k_3y_{t_{n}}^{3/2}\D W_n}{\max \left\{1, \D \left((k_1y_{t_{n}}- k_2y_{t_{n}}^2)\D + k_3y_{t_{n}}^{3/2}\D W_n\right)\right\}}
$$

\begin{figure}[ht]
  \caption{ SD, HMS, and TAMeD method applied to SDE (\ref{eq00001}) with HMS exact solution and parameters $k_1=0.1, k_2=\frac{\lam}{2}(k_3)^2, k_3=\sqrt{0.2}, \lam=700, x_0=1, T=1$ with $17$ digits of accuracy.}
  \centering
    \includegraphics[width=0.5\textwidth]{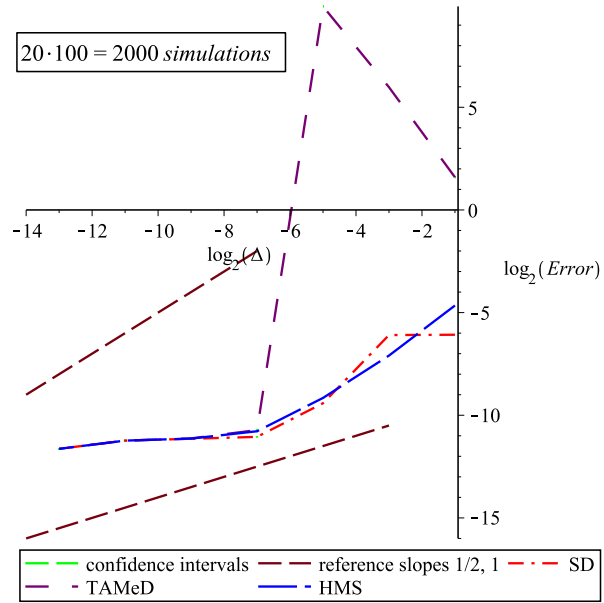}\label{SD_HMS_TDwithHMS_20_100steperror}
\end{figure}

\begin{figure}[ht]
  \caption{ SD, HMS, and TAMeD method applied to SDE (\ref{eq00001}) with SD exact solution and parameters  $k_1=0.1, k_2=\frac{\lam}{2}(k_3)^2, k_3=\sqrt{0.2}, \lam=700, x_0=1, T=1$ with $17$ digits of accuracy.}
  \centering
    \includegraphics[width=0.5\textwidth]{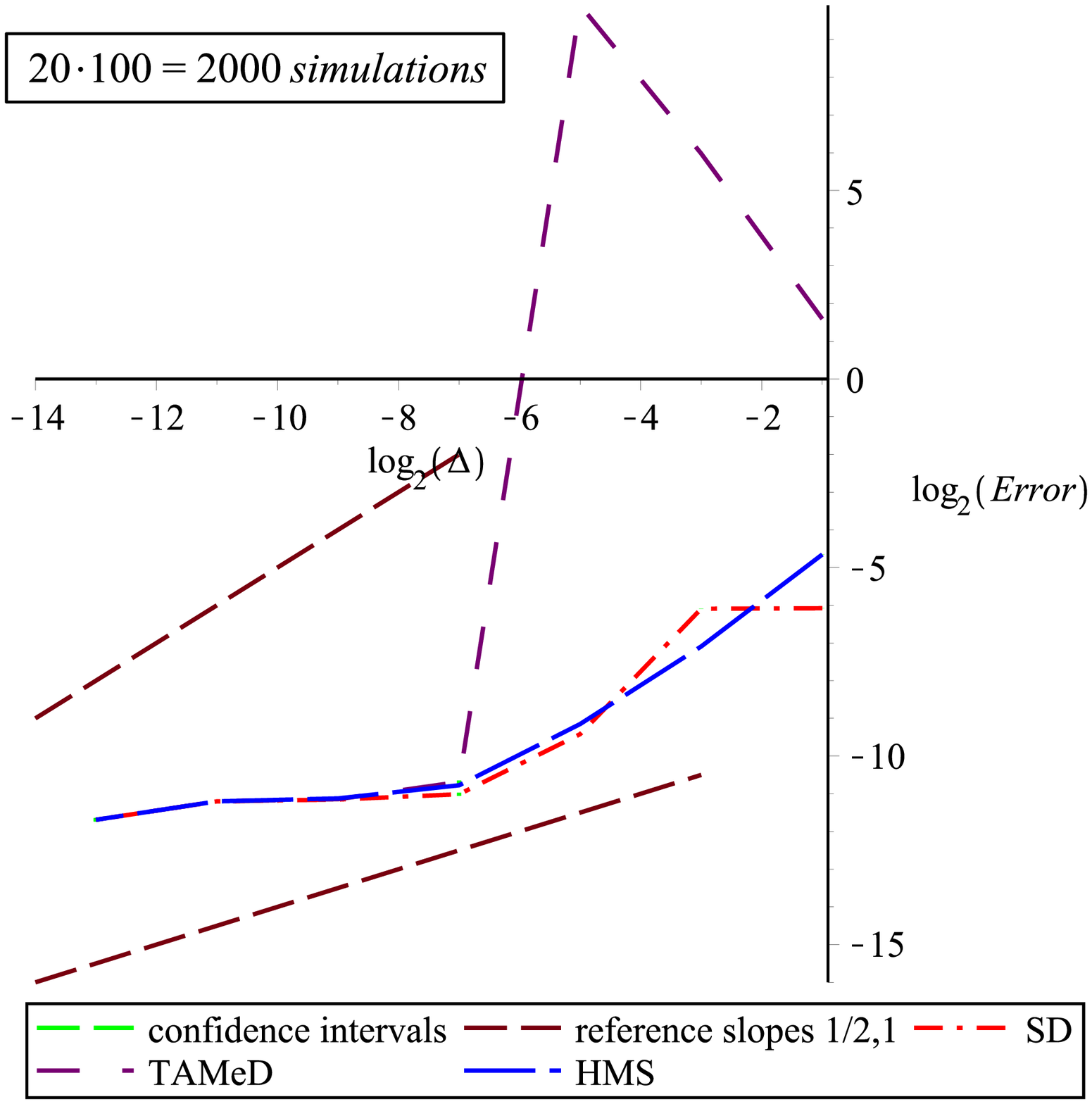}
		\label{SD_HMS_TDwithSD_20_100steperror}
\end{figure}

As a reference solution, we take in the first experiment the value of $y_T^{HMS}$ at $\D=2^{-14},$ as in the numerical experiment in (\cite[Section 4.1]{higham_et_al:2013}), and in the second experiment $y_T^{SD}$ at $\D=2^{-14},$ since we have shown by (\ref{eq00004}) that it strongly converges to the exact solution. We plot in a $\log_2-\log_2$ scale and error bars represent  $90\%$ confidence intervals. The results are shown in Figures \ref{SD_HMS_TDwithHMS_20_100steperror} and \ref{SD_HMS_TDwithSD_20_100steperror} and Tables \ref{tab:withHMS} and \ref{tab:withSD}.

\begin{table}[htbp]  \scriptsize
\centering
		\begin{tabular}{|c|l|l|l|}
		\hline  Step $\D$ & $90\%$ SD-Error & $90\%$ HMS-Error & $90\%$ TAMeD-Error \\	
  \hline $2^{-1}$  & $0.01479749664 \pm 1.584\cdot10^{-5}$  & $0.03968188388\pm 1.610\cdot10^{-5}$ & $3.014797494 \pm 1.584\cdot10^{-5}$\\
      \hline $2^{-3}$  &$0.01464432262  \pm 1.796\cdot10^{-5}$  & $0.007325380970\pm1.810\cdot10^{-5}$ & $ 63.01481485\pm  1.795\cdot10^{-5}$\\
		\hline $2^{-5}$  & $0.001465805974 \pm 1.920\cdot10^{-5}$  & $0.001752988500\pm1.910\cdot10^{-5}$ & $964.1295990\pm 8.992\cdot10^{-3}$\\
		\hline $2^{-7}$  & $0.0004706806728 \pm 1.252\cdot10^{-5} $  & $0.0005690540935\pm1.780\cdot10^{-5}$ & $0.0005921634365\pm 1.677\cdot10^{-5}$\\
		\hline $2^{-9}$  &  $0.0004415939458 \pm 1.311\cdot10^{-5} $ & $0.0004442429779\pm1.385\cdot10^{-5}$ & $0.0004465603424\pm1.319\cdot10^{-5}$\\
		\hline $2^{-11}$ & $0.0004149841292\pm 1.290\cdot10^{-5} $  & $0.0004148866098\pm1.261\cdot10^{-5}$ & $0.0004148921662\pm 1.287\cdot10^{-5}$\\
		\hline $2^{-13}$ &  $0.0003145934380 \pm6.461\cdot10^{-6} $ & $0.0003143683331\pm6.476\cdot10^{-6}$ & $0.0003143198008\pm 6.474\cdot10^{-6}$\\
	\hline	\end{tabular}
	\caption{\small Error and step size of SD,HMS and TAMeD approximation of (\ref{eq00001}) with HMS exact solution with $17$ digits of accuracy.}
	\label{tab:withHMS}
\end{table}

\begin{table}[htbp] \scriptsize
\centering
		\begin{tabular}{|c|l|l|l|}
		\hline  Step $\D$ & $90\%$ SD-Error & $90\%$ HMS-Error & $90\%$ TAMeD-Error \\	
   \hline $2^{-1}$  & $ 0.01478722761\pm1.694\cdot10^{-5}$  & $0.03969436537\pm1.732\cdot10^{-5}$ & $3.014787223\pm1.694\cdot10^{-5}$\\
      \hline $2^{-3}$ & $ 0.01464578986\pm2.049\cdot10^{-5}$  & $0.007323802285\pm2.047\cdot10^{-5}$ & $63.01481630\pm2.050\cdot10^{-5}$\\
		\hline $2^{-5}$  & $0.001460523189\pm1.915\cdot10^{-5}$  & $0.001759496261\pm1.896\cdot10^{-5}$ & $964.1304450\pm6.430\cdot10^{-3}$\\
		\hline $2^{-7}$  & $0.0004839919120\pm1.381\cdot10^{-5}$  & $0.0005708352815\pm1.79\cdot10^{-5}$ & $0.0006062383075\pm1.695\cdot10^{-5}$\\
		\hline $2^{-9}$  &  $ 0.0004393400262\pm1.0818\cdot10^{-5}$ & $0.0004483330032\pm1.108\cdot10^{-5}$ & $0.0004421671951\pm1.0935\cdot10^{-5}$\\
		\hline $2^{-11}$ & $0.0004244777682 \pm1.0218\cdot10^{-5}$  & $0.0004249117572\pm1.018\cdot10^{-5}$ & $0.0004244440682\pm1.021\cdot10^{-5}$\\
		\hline $2^{-13}$ &  $0.0003025212586\pm8.797\cdot10^{-6}$ & $0.0003026818444\pm8.748\cdot10^{-6}$ & $0.0003027212689\pm8.736\cdot10^{-6}$\\	\hline	\end{tabular}
	\caption{\small Error and step size of SD,HMS and TAMeD approximation of (\ref{eq00001}) with SD exact solution with $17$ digits of accuracy.}
	\label{tab:withSD}
\end{table}

The following points of discussion are worth mentioning.
\begin{itemize}
	\item The SD method and the HMS method are very close, with SD performing slightly better, except only for the step size $\D=2^{-3}.$ The same situation appears in both cases, i.e. independently of the choice of the exact solution, which is  a positive feature of SD. 
	\item A linear regression with the method of least squares fit, in the case one considers only the first four points with steps $\D=2^{-1}, 2^{-3}, 2^{-5}, 2^{-7},$ produced values consistent with the strong order of convergence equal to $1$ for both SD and HMS methods, whereas considering all the seven points, values close to $1/2.$ Tables \ref{tab:orderHMS} and \ref{tab:orderSD} present the exact values of order of convergence. We see that the order of convergence of SD for  problem (\ref{eq00001}) is at least $1/2.$  
	\begin{table}[htbp] 
\centering
		\begin{tabular}{|c|l|l|}
		\hline  Number of points & order of SD & order of HMS \\	
    \hline          $4$ &  $0.912$  & $1.022$\\
    \hline          $7$ & $ 0.512$  & $0.557$\\
		\hline
  \end{tabular}
	\caption{\small Order of convergence of SD and HMS approximation of (\ref{eq00001}) with HMS exact solution with $17$ digits of accuracy.}
	\label{tab:orderHMS}
\end{table}
	
	\begin{table}[htbp] 
\centering
		\begin{tabular}{|c|l|l|}
		\hline  Number of points & order of SD & order of HMS \\	
    \hline          $4$ &  $0.906$  & $1.021$\\
    \hline          $7$ & $ 0.514$  & $0.558$\\
		\hline
  \end{tabular}
	\caption{\small Order of convergence of SD and HMS approximation of (\ref{eq00001}) with SD exact solution with $17$ digits of accuracy.}
	\label{tab:orderSD}
\end{table}
	\item The confidence intervals are of such an order that indicates that we donnot need to increase the number of batches $M.$ All the above calculations are made evaluating with $17$ digits. The results of doubling the number of digits to $34$ are shown in the following Tables \ref{tab:withHMS2} and \ref{tab:orderHMS2}, that indicate that there is no significant difference of the situation.
	\begin{table}[htbp]  \footnotesize
\centering
		\begin{tabular}{|c|l|l|l|}
		\hline  Step $\D$ & $90\%$ SD-Error & $90\%$ HMS-Error & $90\%$ TAMeD-Error \\	
  \hline $2^{-1}$  & $0.01480569914 \pm 2.376\cdot10^{-5}$  & $0.03967116854\pm 2.368\cdot10^{-5}$ & $3.014805694\pm 2.376\cdot10^{-5}$\\
      \hline $2^{-3}$  &$0.01462352787\pm1.552\cdot10^{-5}$  & $0.007345838060\pm1.559\cdot10^{-5}$ & $63.01479405\pm1.552\cdot10^{-5}$\\
		\hline $2^{-5}$  & $0.001500299224\pm1.861\cdot10^{-5}$ & $0.001721224885\pm1.838\cdot10^{-5}$ & $964.1293050\pm 9.272\cdot10^{-3}$\\
		\hline $2^{-7}$  & $0.0004733674508\pm1.252\cdot10^{-5} $  & $0.0005777263750\pm1.130\cdot10^{-5}$ & $0.0005898564120\pm 1.339\cdot10^{-5}$\\
		\hline $2^{-9}$  & $0.0004411224169\pm 1.454\cdot10^{-5}$ & $0.0004504228844\pm1.411\cdot10^{-5}$ & $0.0004446873228\pm1.488\cdot10^{-5}$\\
		\hline $2^{-11}$ & $0.0004260658292\pm 1.492\cdot10^{-5}$ & $0.0004255780634\pm1.469\cdot10^{-5}$ & $0.0004258924099\pm1.490\cdot10^{-5}$\\
		\hline $2^{-13}$ &  $0.0003137838988\pm9.182\cdot10^{-6}$ & $0.0003137626410\pm9.144\cdot10^{-6}$ & $0.0003137661728\pm9.139\cdot10^{-6}$\\
	\hline	\end{tabular}
	\caption{\small Error and step size of SD,HMS and TAMeD approximation of (\ref{eq00001}) with HMS exact solution with $34$ digits of accuracy.}
	\label{tab:withHMS2}
\end{table}

		\begin{table}[htbp] 
\centering
		\begin{tabular}{|c|l|l|}
		\hline  Number of points & order of SD & order of HMS \\	
    \hline          $4$ &  $0.909$  & $1.020$\\
    \hline          $7$ & $ 0.512$  & $0.555$\\
		\hline
  \end{tabular}
	\caption{\small Order of convergence of SD and HMS approximation of (\ref{eq00001}) with HMS exact solution with $34$ digits of accuracy.}
	\label{tab:orderHMS2}
\end{table}

		\item For small $\D$ it may happen that the global error will begin to increase as $\D$ is further decreased (\cite[p.97]{kloeden_platen_schurz:2003}). This effect is due to the roundoff error which influences the calculated global error. In practice, that implies the existence of a minimum step size $\D_{\min},$ for each initial value problem, below which the accuracy of the approximations through a specific method cannot be improved. 
	\item Convergence of a numerical scheme does not alone guarantee its practical value (\cite[p.129]{kloeden_platen_schurz:2003}). It may be numerical UNSTABLE. Moreover, in practice, the computer time consumed to provide a desired level of accuracy, is of great importance. As mentioned in Footnote \ref{simu}, we donnot claim that SD method performs well in that aspect, because of the exponential calculations involved. However, it seems that it can reach accuracy up to $4$ digits, as fast as the HMS method.
\item We would like to see how things become, by altering the parameter $\lam.$ SD method, seems to work, with the theoretical proof shown in Section \ref{ssec:3.1}, when $\lam$ is over $7.$ What happens below that range? HMS method works for $\lam$ over $1/2.$ Moreover, as noted in Remark \ref{r00001}(iv), our method can cover more general cases, in contrast to HMS, by introducing the function $\phi(\cdot)$ in the diffusion part, or/and by assuming random coefficients $k_1(\cdot),k_2(\cdot),k_3(\cdot).$ 

In the following Figure \ref{SD_HMS_TDwithHMS_20_100steperror_C=10} we present the situation when we change the parameters of SDE (\ref{eq00001}) in such a way that we are closer to the theoretical acceptable range ( by lowering $\lam$ to $70$). 
\begin{figure}[ht]
  \caption{ SD, HMS, and TAMeD method applied to SDE (\ref{eq00001}) with HMS exact solution and parameters $k_1=0.1, k_2=\frac{\lam}{2}(k_3)^2, k_3=\sqrt{0.2}, \lam=70, x_0=1, T=1$ with $17$ digits of accuracy.}
  \centering
    \includegraphics[width=0.5\textwidth]{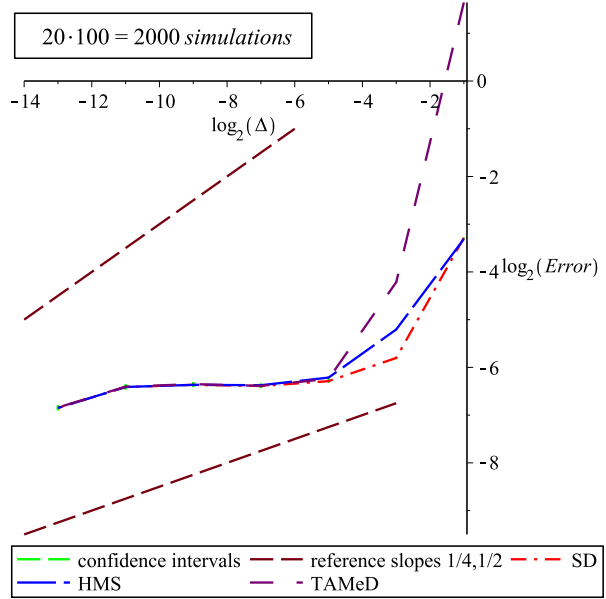}\label{SD_HMS_TDwithHMS_20_100steperror_C=10}
\end{figure}
The rate of convergence drops to a half for both SD and HMS method and TAMeD seems to perform better than before. To be more precise we present in the table \ref{tab:orderHMS10} the exact numbers.
\begin{table}[htbp] 
\centering
		\begin{tabular}{|c|l|l|}
		\hline  Number of points & order of SD & order of HMS \\	
    \hline          $4$ &  $0.490$  & $0.510$\\
    \hline          $7$ & $ 0.214$  & $0.235$\\
		\hline
  \end{tabular}
	\caption{\small Order of convergence of SD and HMS approximation of (\ref{eq00001}) with HMS exact solution with $17$ digits of accuracy when $\lam=70.$}
	\label{tab:orderHMS10}
\end{table}

In Figure \ref{SD_HMS_TDwithHMS_20_100steperror_C=1} we present the case with $\lam=7.$
\begin{figure}[ht]
  \caption{ SD, HMS, and TAMeD method applied to SDE (\ref{eq00001}) with HMS exact solution and parameters $k_1=0.1, k_2=\frac{\lam}{2}(k_3)^2, k_3=\sqrt{0.2}, \lam=7, x_0=1, T=1$ with $17$ digits of accuracy.}
  \centering
    \includegraphics[width=0.5\textwidth]{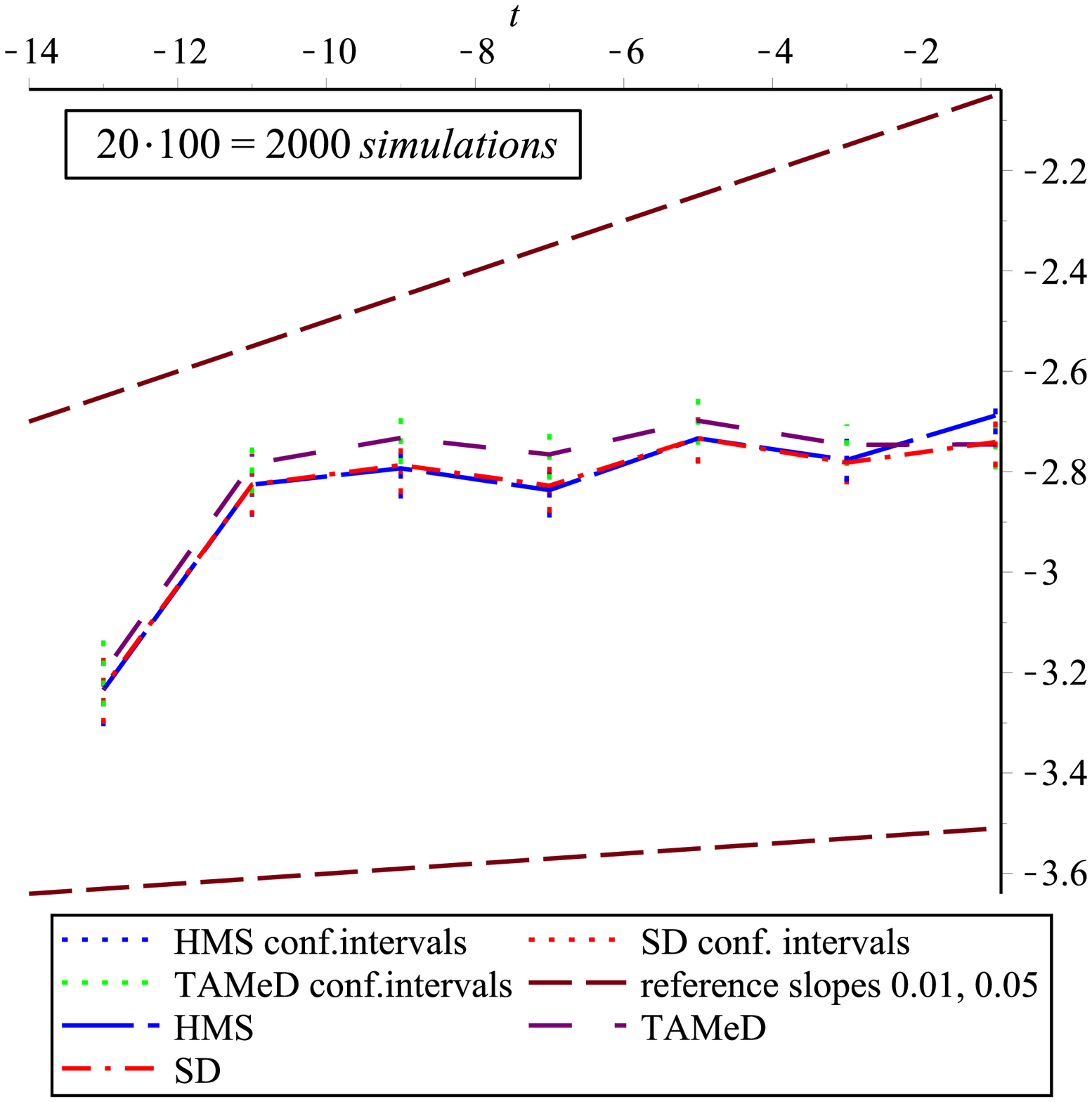}\label{SD_HMS_TDwithHMS_20_100steperror_C=1}
\end{figure}
The rate of convergence drops dramatically for all methods. Moreover the TAMeD performs even better, close to SD and HMS. To be more precise we present in the table \ref{tab:orderHMS1} the exact numbers.
\begin{table}[htbp] 
\centering
		\begin{tabular}{|c|l|l|l|}
		\hline  Number of points & order of SD & order of HMS  & order of TAMeD\\	
    \hline          $7$ & $ 0.029$  & $0.032$  & $0.026$\\
		\hline
  \end{tabular}
	\caption{\small Order of convergence of SD and HMS approximation of (\ref{eq00001}) with HMS exact solution with $17$ digits of accuracy when $\lam=7.$}
	\label{tab:orderHMS1}
\end{table}

\item Regarding the TAMeD method, a major drawback is that it does not preserve positivity. However,  we remark that even though the errors of the TAMeD approximation are quite big, for big step sizes,\footnote{In the plots these errors donnot seem so big, because of the $\log_2-$scale. The tables though show this anomaly.} all methods behave quite close for small $\D$'s and even closer for bigger $\D$ as we lower the parameter $\lam$ close to its critical value.
\end{itemize}



\appendix
\section{Existence and uniqueness of $y_t^{SD}$ for Heston $3/2-$model}\label{ap1}
\subsection{Uniqueness of solution of $y_t^{SD}$}\label{ap1.1}
\setcounter{equation}{0}

Let $y_t, \hat{y}_t$ be two solutions of SDE (\ref{eq00003}) with same initial condition, i.e. with $y_0=\hat{y}_0.$ By Lemma \ref{l3} they both belong to the space $\bbm^2([0,T];\bbR)$ of measurable $\{\bbf_t\}$ adapted processes $z$ such that
$$
\bfE \int_0^T|z_s|^2ds<\infty.
$$
Set the stopping times $\theta^i_R=\inf\{t\in [t_{i-1},t_{i}]: |y_t|>R\}$ and $\hat{\theta}^i_R=\inf\{t\in [t_{i-1},t_{i}]: |\hat{y}_t|>R\}$ for some $R>0$ big enough and consider the stopping times $\tau^i_R=\theta_R^i\wedge\hat{\theta}_R^i,$ for $i=1,...,N.$
Take $t\in[0,t_1]$ and $e_{t\wedge\tau_R^1}:=y_{t\wedge\tau_R^1}-\hat{y}_{t\wedge\tau_R^1}.$ It holds that
\beao
|e_{t\wedge\tau_R^1}|^2&=&\!\!\!\!\!\left| \int_{0}^{t\wedge\tau_R^1}\left(f(\hat{s},s,y_{\hat{s}},y_s)-f(\hat{s},s,\hat{y}_{\hat{s}},\hat{y}_s)\right)ds + \int_{0}^{t\wedge\tau_R^1}\left(g(\hat{s},s,y_{\hat{s}},y_s) - g(\hat{s},s,\hat{y}_{\hat{s}},\hat{y}_s)\right) dW_s\right|^2\\
&\leq&2t\int_{0}^{t_1\wedge\tau_R^1}\Big|f(\hat{s},s,y_{\hat{s}},y_s)-f(\hat{s},s,\hat{y}_{\hat{s}},\hat{y}_s)\Big|^2ds + 2|M_{t}|^2\\
&\leq& 2t\int_{0}^{t\wedge\tau_R^1}4C_R^2 \left(|y_{\hat{s}}-\hat{y}_{\hat{s}}|^2 + |y_s-\hat{y}_s|^2 + |y_{\hat{s}}-\hat{y}_{\hat{s}}|^{2\rho}\right) ds + 2 |M_{t}|^2\\
&\leq& 8tC_R^2\int_{0}^{t}|e_{s\wedge\tau_R^1}|^2ds + 2 |M_{t}|^2,
\eeao
where in the second step Cauchy-Schwarz inequality, in the third step the elementary inequality $(\sum_{i=1}^3 a_i)^2 \leq 4\sum_{i=1}^3 a_i^2,$ for the appropriate $a_i$'s and Assumption A for $f,$ in the last step the fact that $\hat{s}=0,$ when $s\in[0,t_1]$ and the equality in the initial conditions $y_0=\hat{y}_0$ and
$$
M_{t}:=\int_{0}^{t\wedge\tau_R^1}\left( g(\hat{s},s,y_{\hat{s}},y_s)-g(\hat{s},s,\hat{y}_{\hat{s}},\hat{y}_s) \right)dW_s.
$$
Taking the supremum over all $t\in[0,t_1]$ and then expectations we have
\beam
\nonumber
&&\bfE\sup_{0\leq t\leq t_1}|e_{t\wedge\tau_R^1}|^2
\leq 8t C_R^2 \bfE\sup_{0\leq t\leq t_1}\left(\int_{0}^{t\wedge\tau_R^1}|y_s-\hat{y}_s|^2ds\right)+ 2\bfE\sup_{0\leq t\leq t_1} |M_t|^2\\
\label{eq19}&\leq& 8t_1C_R^2 \int_{0}^{t_1}\bfE\sup_{0\leq l\leq s}|e_{l\wedge\tau_R^1}|^2ds + 2\bfE |M_{t_1}|^2,
\eeam
where we have used Doob's maximal inequality with $p=2,$ since  $M_t$ is an $\bbR-$valued martingale that belongs to $\bbl^2.$ Moreover, we have that
\beao
&&\bfE |M_{t_1}|^2:=\bfE\left|\int_{0}^{t_1\wedge\tau_R^1} \left( g(\hat{s},s,y_{\hat{s}},y_s)-g(\hat{s},s,\hat{y}_{\hat{s}},\hat{y}_s) \right) dW_s\right|^2\\
&=&
 \bfE\left(\int_{0}^{t_1\wedge\tau_R^1} \left( g(\hat{s},s,y_{\hat{s}},y_s)-g(\hat{s},s,\hat{y}_{\hat{s}},\hat{y}_s) \right)^2ds\right)\\
&\leq& 4C_R^2\bfE\left(\int_{0}^{t_1\wedge\tau_R^1} \left( |y_{0}-\hat{y}_{0}|^2 + |y_s-\hat{y}_s|^2 + |y_{0}-\hat{y}_{0}|\right) ds\right)\\
&\leq& 4C_R^2\int_{0}^{t_1\wedge\tau_R^1} \bfE|y_s-\hat{y}_{s}|^2ds \leq 4C_R^2 \int_{0}^{t_1}\bfE\sup_{0\leq l\leq s}|e_{l\wedge\tau_R^1}|^2ds,
\eeao
where we have used Assumption A for $g,$ thus relation (\ref{eq19}) becomes
$$
\bfE\sup_{0\leq t\leq t_1}|e_{t\wedge\tau_R^1}|^2 \leq (8t_1C_R^2 + 4C_R^2)\int_{0}^{t_1}\bfE\sup_{0\leq l\leq s}|e_{l\wedge\tau_R^1}|^2ds,
$$
which by use of Gronwall's inequality gives
\beqq\label{eq20}
\bfE\sup_{0\leq t\leq t_1}|e_{t\wedge\tau_R^1}|^2=0.
\eeqq
Following the same arguments we can show that
$$
\bfE\sup_{0\leq t\leq t_1}|e_{t\wedge\tau_R^i}|^2=0,
$$
for every integer $1\leq i\leq N.$\footnote{For $i=2$ just use the same ideas as for $i=1$  and the other cases follow exactly the same way using in every step the result of the previous step.}
Thus, if we drop the index $i$ from the stopping times with the meaning that $\theta_R=\inf\{t\in [0,T]: |y_t|>R\}$ and $\hat{\theta}_R=\inf\{t\in [0,T]: |\hat{y}_t|>R\}$ for some $R>0$ big enough and consider the stopping time $\tau_R=\theta_R\wedge\hat{\theta}_R,$  we have that
$$
\bfE\sup_{0\leq t\leq T}|e_{t\wedge\tau_R}|^2\leq \sum_{i=1}^N \bfE\sup_{t_{i-1}\leq t\leq t_i}|e_{t\wedge\tau_R^i}|^2=0.
$$
Hence, $y_t=\hat{y}_t$ for all $0\leq t\leq T$ a.s. which proves that the solution of SDE (\ref{eq00003}), and in general of  SDE (\ref{eq1.1}) when it exists, is unique.

\subsection{Existence of solution of $y_t^{SD}$}\label{ap1.2}

We will show the existence of the solution of SDE (\ref{eq00002}) for $n=0$ and the same procedure can be followed to show  the existence of the solution of SDE (\ref{eq00002}) for every integer $n=1,..,N-1,$ i.e. the existence of the solution of SDE (\ref{eq00003}). Application of Ito's formula to $\ln y_t,$ for $0\leq t\leq t_1$ implies
\beao
&&\ln y_t= \ln y_0 + \int_{0}^{t} \frac{1}{y_s}(k_1(s) - k_2(s)y_{0})y_sds + \frac{1}{2}\int_{0}^{t}\left(-\frac{1}{y_s^2}\right) k_3^2(s)y_{0}y_s^2 ds\\
&&+ \int_{0}^{t}\frac{1}{y_s} k_3(s)y_{0}y_s dW_s\\
&=&\ln y_0 + \int_{0}^{t} \left(k_1(s) - k_2(s)y_{0} -\frac{k_3^2(s)}{2}\sqrt{y_{0}}\right) ds + \int_{0}^{t} k_3(s)\sqrt{y_{0}}dW_s.
\eeao
Now take the exponential of both sides of (\ref{eq14}) with $\hat{s}=0$ in the case $0\leq t\leq t_1$ to verify that (\ref{eq00005}) is indeed a solution of SDE (\ref{eq00002}) for $n=0.$

\end{document}